\documentclass[11pt]{article}
\usepackage{amssymb}

\usepackage{amsmath}
\usepackage{theorem}
\usepackage{here}
\usepackage[pdftex]{graphicx}

\newtheorem{thm}{Theorem}[section]
\newtheorem{lem}{Lemma}[section]
\newtheorem{cor}{Corollary}[section]
\newtheorem{prp}{Proposition}[section]
\theorembodyfont{\rmfamily}

\makeatletter

\@addtoreset{equation}{section}
\makeatother
\def\Yt{{\widetilde Y}}

\usepackage{comment} %
\usepackage{bm}

\def\ta{{\tau}}

\def\Xt{{\widetilde X}}
\def\Yt{{\widetilde Y}}

\def\pd{\partial}

\def\al{{\alpha}}
\def\be{{\beta}}
\def\ga{{\gamma}}
\def\de{{\delta}}

\def\la{{\lambda}}
\def\si{{\sigma}}

\def\ze{{\zeta}}

\def\lah{{\hat \la}}

\def\De{{\Delta}}

\def\fh{{\hat f}}
\def\gh{{\hat g}}

\def\ph{{\hat p}}

\def\[{{\text{\boldmath $[$}}}
\def\]{{\text{\boldmath $]$}}}

\def\/{{\Bigr/\!\!}}

\def\1r{{\rm (1)}}
\def\2r{{\rm (2)}}
\def\3r{{\rm (3)}}
\def\4r{{\rm (4)}}
\def\5r{{\rm (5)}}

\def\non{{\nonumber}}

\begin{document}
\title{Bayesian Point Estimation and Predictive Density Estimation for the Binomial Distribution with a Restricted Probability Parameter}
\author{
Yasuyuki Hamura\footnote{Graduate School of Economics, University of Tokyo, 
7-3-1 Hongo, Bunkyo-ku, Tokyo 113-0033, JAPAN. 
JSPS Research Fellow.\newline{
E-Mail: yasu.stat@gmail.com}} \
}
\maketitle
\begin{abstract}
In this paper, we consider Bayesian point estimation and predictive density estimation in the binomial case. 
After presenting preliminary results on these problems, we compare the risk functions of the Bayes estimators based on the truncated and untruncated beta priors and obtain dominance conditions when the probability parameter is less than or equal to a known constant. 
The case where there are both a lower bound restriction and an upper bound restriction is also treated. 
Then our problems are shown to be related to similar problems in the Poisson case. 
Finally, numerical studies are presented. 

\par\vspace{4mm}
\end{abstract}

\section{Introduction}
\label{sec:introduction}
Estimation of a restricted parameter has been widely studied in the literature. 
For example, Casella and Strawderman (1981), Marchand and Perron (2001), and Hartigan (2004) considered the normal case. 
Singh, Gupta and Nisra (1993), Parsian and Sanjari Farsipour (1997), and Tripathi, Kumar and Petropoulos (2014) considered the exponential case. 
The binomial case was treated by Perron (2003). 
Recently, Marchand, Rancourt and Strawderman (2021) treated one-parameter exponential families. 
See, for example, van Eeden (2006) for more details. 

Similar results in the context of predictive density estimation under the Kullback-Leibler (KL) divergence have also been obtained. 
Fourdrinier et al. (2011), Hamura and Kubokawa (2019), and Hamura and Kubokawa (2020) treated the normal, negative binomial, and Poisson cases by using identies which relate Bayesian predictive density estimation to Bayesian point estimation. 
L'Moudden et al. (2017) and Hamura and Kubokawa (2021) considered the gamma and exponential cases directly. 
General results for location and/or scale families were obtained by Kubokawa et al. (2013). 

In the binomial case, although a simple identity relating Bayesian predictive density estimation to Bayesian point estimation is available, the former has not been fully considered when the parameter space is restricted. 
Some distinctive properties of the binomial distribution is as follows: 
\begin{itemize}
\item
The Jeffreys prior is proper and hence the associated Bayesian procedures are admissible when the parameter space is unrestricted. 
\item
The entropy loss in the binomial case is neither the squared error loss nor the entropy loss in the gamma case nor their weighted versions. 
\item
The sample space is a finite set. 
\end{itemize}
In this paper, we show that domination results hold in the binomial case also. 
In addition, we discuss similarities and dissimilarities between the binomial and Poisson cases. 

In the $\text{binomial} (n, p)$ case, there are different default choices of a prior distribution $\pi (p)$ for the probability parameter $p$. 
For example, the Jeffreys prior is $\pi (p) \propto p^{- 1 / 2} (1 - p)^{- 1 / 2}$. 
Tuyl, Gerlachy and Mengersen (2009) argued for the uniform prior $\pi (p) \propto 1$. 
Komaki (2012) showed that the Bayesian predictive density based on the prior $\pi (p) \propto p^{1 / \sqrt{6}} (1 - p)^{1 / \sqrt{6}}$ is asymptotically minimax. 
The prior $\pi (p) \propto p^{\sqrt{n} / 2 - 1} (1 - p)^{\sqrt{n} / 2 - 1}$ yields the unique minimax estimator under the squared error loss (see, for example, Lehmann and Casella (1998)). 
The improper prior $\pi (p) \propto p^{- 1} (1 - p)^{- 1}$ cannot be used since it can lead to an improper posterior. 
Therefore, the general conjugate beta prior $\pi (p) \propto \pi _{a, b} (p) = p^{a - 1} (1 - p)^{b - 1}$ for $a, b > 0$ is treated as a benchmark in this paper. 

In Section \ref{sec:domination}, we consider the situation where for some known $0 < \overline{p} < 1$, the probability parameter $p$ is restricted to the interval $(0, \overline{p} ]$. %
In this situation, it seems natural to use the truncated beta prior $\pi _{a, b, \overline{p}} (p) = \pi _{a, b} (p) 1_{(0, \overline{p} ]} (p)$ instead of $\pi _{a, b} (p)$. 
However, as in the cases of other distributions, it is not necessarily easy to determine sufficient conditions under which the Bayesian procedures based on $\pi _{a, b, \overline{p}} (p)$ dominate those based on $\pi _{a, b} (p)$. 
In fact, later we show that the result does not always hold, no matter what the values of $a, b > 0$ are. 
Thus, obtaining sufficient conditions is important. 

The remainder of the paper is organized as follows. 
In Section \ref{sec:preliminaries}, preliminary results on Bayesian point estimation and predictive density estimation in the binomial case are presented. 
In Section \ref{sec:domination}, we obtain conditions for the Bayes estimator based on $\pi _{a, b, \overline{p}} (p)$ to dominate that based on $\pi _{a, b} (p)$ when $p \in (0, \overline{p} ]$. 
In Section \ref{sec:both}, we consider the case where there are both a lower bound restriction and an upper bound restriction. 
In Section \ref{sec:relation_Po}, we show that our problems are related to similar problems in the Poisson case. 
Numerical studies are presented in Section \ref{sec:sim}. 
All lemmas are proved in the Appendix.

\def\pit{{\tilde \pi}}

\section{Preliminaries}
\label{sec:preliminaries}
In this section, preliminary results on Bayesian point estimation and predictive density estimation in the binomial case are presented. 
We consider the KL divergence and the entropy loss derived from the KL divergence. 
For $n \in \mathbb{N}$ and $p \in (0, 1)$, let 
\begin{align}
f(x | n, p) &= \binom{n}{x} p^x (1 - p)^{n - x} \text{,} \quad x = 0, 1, \dots , n \text{,} \non 
\end{align}
denote the density of the binomial distribution ${\rm{Bin}} (n, p)$.

\subsection{Plug-in predictive density estimation}
\label{subsec:plug-in}
For known $n, l \in \mathbb{N}$ and unknown $p \in (0, 1)$, let $X \sim {\rm{Bin}} (n, p)$ and $Y \sim {\rm{Bin}} (l, p)$ be independent current and future binomial variables. 
Consider the problem of estimating the predictive density $f( \cdot | l, p)$ of $Y$ on the basis of $X$ under the KL divergence. 
Then the risk function of a predictive density estimator $\fh ( \cdot ; X)$ is 
\begin{align}
R_{l, n} (p, \fh ) &= E_{p}^{(Y, X)} \Big[ \log {f(Y | l, p) \over \fh (Y ; X)} \Big] \text{.} \non 
\end{align}
For an estimator $\de (X) \in (0, 1)$ for $p$, we write $\fh _{l}^{( \de )} ( \cdot ; X)$ for the plug-in density estimator $f( \cdot | l, \de (X))$. 
As shown by Robert (1996), the risk function of $\fh _{l}^{( \de )} ( \cdot ; X)$ is simply $l$ times the risk function of $\de (X)$ under the entropy loss; that is, 
\begin{align}
R_{l, n} (p, \fh _{l}^{( \de )} ) &= %
E_{p}^{X} \Big[ l \Big\{ p \log {p \over \de (X)} + (1 - p) \log {1 - p \over 1 - \de (X)} \Big\} \Big] \non \\
&= E_{p}^{X} [ l L( \de (X), p) ] = l R_n (p, \de ) = l R_{1,n} (p, \fh _{1}^{( \de )} ) \text{,} \non 
\end{align}
where 
\begin{align}
L(d, p) &= p \log {p \over d} + (1 - p) \log {1 - p \over 1 - d} \text{,} \quad d \in (0, 1) \text{,} \label{eq:loss_entropy} 
\end{align}
is the entropy loss in the binomial case and where $R_n (p, \de )$ is the risk function of $\de (X)$ under the entropy loss. 
In contrast to the Poisson case, both terms in the above expression are nonlinear in $d$ and $p$.

\subsection{Bayesian point estimation}
\label{subsec:estimation}
The usual estimator of $p$ under the squared error loss is the maximum likelihood estimator $X / n$. 
However, since it can take on the values $0, 1$ with positive probability, we cannot use it in constructing a plug-in density estimator (as noted in Section 4 of Aitchison (1975)) or in estimating $p$ under the entropy loss. 
A convenient way to obtain an estimator taking on values in $(0, 1)$ is through the use of a prior for $p$. 
The Bayes estimator with respect to a prior $\pi (p)$ and the entropy loss is the posterior mean 
\begin{align}
\ph _{n}^{( \pi )} (X) &= E_{\pi }^{p | X} [ p | X ] = \frac{ \int_{0}^{1} p^{1 + X} (1 - p)^{n - X} \pi (p) dp }{ \int_{0}^{1} p^{X} (1 - p)^{n - X} \pi (p) dp } \text{.} \non 
\end{align}
It is in $(0, 1)$ with probability one. 
Therefore, its risk function is defined under the entropy loss.

\subsection{Bayesian predictive density estimation}
\label{subsec:prediction}
While $\ph _{n}^{( \pi )} (X)$ could be used to construct the plug-in density estimator $\fh _{l}^{( \ph _{n}^{( \pi )} )} ( \cdot ; X) = f( \cdot | l, \ph _{n}^{( \pi )} (X))$, a more natural density estimator based on the prior $\pi (p)$ is the Bayesian predictive density obtained by calculating the posterior mean of $f(y | l, p)$ for each $y \in \{ 0, 1, \dots , l \} $. 
More specifically, the Bayesian predictive density, denoted by $\fh _{l, n}^{( \pi )} ( \cdot ; X)$, is given by 
\begin{align}
\fh _{l, n}^{( \pi )} (y; X) &= E_{\pi }^{p | X} [ f(y | l, p) | X ] = \binom{l}{y} \frac{ \int_{0}^{1} p^{y + X} (1 - p)^{l - y + n - X} \pi (p) dp }{ \int_{0}^{1} p^{X} (1 - p)^{n - X} \pi (p) dp } \text{.} \non 
\end{align}
Aitchison (1975) showed that it is the Bayes solution with respect to the prior $\pi (p)$ under the KL divergence. 
Its functional form can be different from that of the binomial density function.

\subsection{A connection formula}
\label{subsec:connection}
The risk function of the Bayesian predictive density $\fh _{l, n}^{( \pi )} ( \cdot ; X)$ is 
\begin{align}
R_{l, n} (p, \fh _{l, n}^{( \pi )} ) &= %
E_{p}^{(Y, X)} \Big[ l \{ p \log p + (1 - p) \log (1 - p) \} - \log \frac{ \int_{0}^{1} \xi ^{Y + X} (1 - \xi )^{l + n - (Y + X)} \pi ( \xi ) d\xi }{ \int_{0}^{1} \xi ^{X} (1 - \xi )^{n - X} \pi ( \xi ) d\xi } \Big] \text{.} \non 
\end{align}
Now let $X_i \sim {\rm{Bin}} (n + i, p)$ for $i = 0, 1, \dots , l$ and let $Z \sim {\rm{Bin}} (1, p)$ be an independent Bernoulli variable. 
Then, since $X \stackrel{{\rm{d}}}{=} X_0$, $X + Y \stackrel{{\rm{d}}}{=} X_l$, and $X_{i + 1} \stackrel{{\rm{d}}}{=} X_i + Z$ for $i = 0, \dots , l - 1$, we have 
\begin{align}
&E_{p}^{(Y, X)} \Big[ \log \frac{ \int_{0}^{1} \xi ^{Y + X} (1 - \xi )^{l + n - (Y + X)} \pi ( \xi ) d\xi }{ \int_{0}^{1} \xi ^{X} (1 - \xi )^{n - X} \pi ( \xi ) d\xi } \Big] \non \\
&= \sum_{i = 0}^{l - 1} E_{p}^{( X_{i + 1} , X_i )} \Big[ \log \frac{ \int_{0}^{1} \xi ^{X_{i + 1}} (1 - \xi )^{n + i + 1 - X_{i + 1}} \pi ( \xi ) d\xi }{ \int_{0}^{1} \xi ^{X_i} (1 - \xi )^{n + i - X_i} \pi ( \xi ) d\xi } \Big] \non \\
&= \sum_{i = 0}^{l - 1} E_{p}^{(Z, X_i )} \Big[ \log \frac{ \int_{0}^{1} \xi ^{Z + X_i} (1 - \xi )^{1 + n + i - (Z + X_i )} \pi ( \xi ) d\xi }{ \int_{0}^{1} \xi ^{X_i} (1 - \xi )^{n + i - X_i} \pi ( \xi ) d\xi } \Big] \non \\
&= \sum_{i = 0}^{l - 1} E_{p}^{X_i} \Big[ \sum_{z = 0}^{1} f(z | 1, p) \log \frac{ \int_{0}^{1} \xi ^{z + X_i} (1 - \xi )^{1 + n + i - (z + X_i )} \pi ( \xi ) d\xi }{ \int_{0}^{1} \xi ^{X_i} (1 - \xi )^{n + i - X_i} \pi ( \xi ) d\xi } \Big] \text{.} \non 
\end{align}
This leads to the following fact. 
\begin{thm}
\label{thm:connection} 
The risk function of the Bayesian predictive density $\fh _{l, n}^{( \pi )} ( \cdot ; X)$ can be expressed using the risk functions of the Bayes estimators $\ph _{n}^{( \pi )} ( X_0 ), \dots , \ph _{n + l - 1}^{( \pi )} ( X_{l - 1} )$ as 
\begin{align}
R_{l, n} (p, \fh _{l, n}^{( \pi )} ) &= \sum_{i = 0}^{l - 1} R_{n + i} (p, \ph _{n + i}^{( \pi )} ) \text{.} \non 
\end{align}
Equivalently, 
\begin{align}
R_{l, n} (p, \fh _{l, n}^{( \pi )} ) &= \sum_{i = 0}^{l - 1} R_{1, n + i} (p, \fh _{1, n + i}^{( \pi )} ) \text{,} \non 
\end{align}
where, for each $i = 0, \dots , l - 1$, $R_{1, n + i} (p, \fh _{1, n + i}^{( \pi )} )$ is the risk function of the one-step-ahead Bayesian predictive density $\fh _{1, n + i}^{( \pi )} ( \cdot ; X_i )$ at time $n + i$. 
\end{thm}

The first equation in the above theorem links Bayesian predictive density estimation to Bayesian point estimation in the binomial case. 
The relation $R_{n} (p, \ph _{n}^{( \pi )} ) = R_{1, n} (p, \fh _{1, n}^{( \pi )} )$ is similar to the relation $R_n (p, \de ) = R_{1,n} (p, \fh _{1}^{( \de )} )$.

\section{The Upper-Bound-Restriction Case}
\label{sec:domination}
In this section, we assume that $p \in (0, \overline{p} ]$ for $0 < \overline{p} < 1$. 
We focus on the conjugate beta prior $\pi _{a, b} (p) = p^{a - 1} (1 - p)^{b - 1}$ and compare it with the truncated version $\pi _{a, b, \overline{p}} (p) = \pi _{a, b} (p) 1_{(0, \overline{p} ]} (p)$, where $a, b > 0$. 
Since, by Theorem \ref{thm:connection}, the risk function of an arbitrary Bayesian predictive density can be expressed as the sum of the risk functions of the corresponding Bayes estimators, we only consider the point estimation of $p$ on the basis of $X \sim {\rm{Bin}} (n, p)$. 
For notational simplicity, we write $\ph ^{( \pi )} = \ph _{n}^{( \pi )} (X)$ for a prior $\pi (p)$ for $p$. 
We write $R(p, \de ) = R_n (p, \de ) = E_{p}^{X} [ L(\de (X), p) ]$ for an estimator $\de (X)$ of $p$.

\subsection{Bayes estimators}
\label{subsec:estimators}
In this section, we derive the Bayes estimators of $p$ based on the priors $\pi _{a, b} (p)$ and $\pi _{a, b, \overline{p}} (p)$. 
First, since $p | X \sim {\rm{Beta}} (X + a, n - X + b)$ under $\pi _{a, b} (p)$, we have $\ph ^{( \pi _{a, b} )} = (X + a) / (n + a + b)$. 
Next, we consider $\ph ^{( \pi _{a, b, \overline{p}} )}$. 
Let 
\begin{align}
I( \al , \ga , \overline{p} ) &= \int_{0}^{1} {t^{\al - 1} \over \{ 1 - \overline{p} (1 - t) \} ^{\ga }} dt \non 
\end{align}
for $\ga > \al > 0$. 

\begin{prp}
\label{prp:estimator_truncated} 
The Bayes estimator $\ph ^{( \pi _{a, b} )}$ is given by 
\begin{align}
\ph ^{( \pi _{a, b, \overline{p}} )} = \ph ^{( \pi _{a, b} )} - {1 \over (n + a + b) I(X + a, n + a + b, \overline{p} )} \text{.} \non 
\end{align}
\end{prp}

\noindent
{\bf Proof%
.} \ \ We have 
\begin{align}
\ph ^{( \pi _{a, b, \overline{p}} )} &= %
\frac{ \int_{0}^{\overline{p}} p^{X + a} (1 - p)^{n - X + b - 1} dp }{ \int_{0}^{\overline{p}} p^{X + a - 1} (1 - p)^{n - X + b - 1} dp } \text{.} \non 
\end{align}
By integration by parts, 
\begin{align}
&- (n - X + b) \int_{0}^{\overline{p}} p^{X + a} (1 - p)^{n - X + b - 1} dp \non \\
&= \Big[ p^{X + a} (1 - p)^{n - X + b} \Big] _{0}^{\overline{p}} - (X + a) \int_{0}^{\overline{p}} p^{X + a - 1} (1 - p)^{n - X + b} dp \non \\
&= \overline{p} ^{X + a} (1 - \overline{p} )^{n - X + b} - (X + a) \Big\{ \int_{0}^{\overline{p}} p^{X + a - 1} (1 - p)^{n - X + b - 1} dp - \int_{0}^{\overline{p}} p^{X + a} (1 - p)^{n - X + b - 1} dp \Big\} \text{.} \non 
\end{align}
Therefore, 
\begin{align}
\ph ^{( \pi _{a, b, \overline{p}} )} &= {X + a \over n + a + b} - {\overline{p} ^{X + a} (1 - \overline{p} )^{n - X + b} \over n + a + b} / \int_{0}^{\overline{p}} p^{X + a - 1} (1 - p)^{n - X + b - 1} dp \text{.} \non 
\end{align}
Since 
\begin{align}
\int_{0}^{\overline{p}} p^{X + a - 1} (1 - p)^{n - X + b - 1} dp &= \int_{0}^{\overline{p} / (1 - \overline{p})} {u^{X + a - 1} \over (1 + u)^{n + a + b}} du \non \\
&= \int_{0}^{1} \Big( {\overline{p} \over 1 - \overline{p}} \Big) ^{X + a} {t^{X + a - 1} \over [1 + \{ \overline{p} / (1 - \overline{p} ) \} t]^{n + a + b}} dt \text{,} \non 
\end{align}
it follows that 
\begin{align}
\ph ^{( \pi _{a, b, \overline{p}} )} &= {X + a \over n + a + b} - {1 \over n + a + b} / \int_{0}^{1} {t^{X + a - 1} \over \{ 1 - \overline{p} (1 - t) \} ^{n + a + b}} dt \text{,} \non 
\end{align}
which is the desired result. 
\hfill$\Box$

\bigskip

We have $\ph ^{( \pi _{a, b, \overline{p}} )} \in (0, \overline{p} )$ by definition, whereas $\ph ^{( \pi _{a, b} )} %
\in [a / (n + a + b), (n + a) / (n + a + b)]$. 
Also, $\lim_{\overline{p} \to 0} \ph ^{( \pi _{a, b, \overline{p}} )} = 0$ and $\lim_{\overline{p} \to 1} \ph ^{( \pi _{a, b, \overline{p}} )} = \ph ^{( \pi _{a, b} )}$. 
It follows from Proposition \ref{prp:estimator_truncated} that 
\begin{align}
{\ph ^{( \pi _{a, b, \overline{p}} )} \over \ph ^{( \pi _{a, b} )}} &= 1 - {1 \over (X + a) I(X + a, n + a + b, \overline{p} )} \label{eq:ratio_p} 
\end{align}
and that 
\begin{align}
{1 - \ph ^{( \pi _{a, b, \overline{p}} )} \over 1 - \ph ^{( \pi _{a, b} )}} &= 1 + {1 \over (n - X + b) I(X + a, n + a + b, \overline{p} )} \text{.} \label{eq:ratio_q} 
\end{align}
Since the integral $I(X + a, n + a + b)$ can be expressed using the incomplete beta function, numerically calculating the estimator $\ph ^{( \pi _{a, b, \overline{p}} )}$ is relatively easy.

\subsection{Sufficient conditions for domination}
\label{subsec:sufficient}
In this section, we compare the risk functions of $\ph ^{( \pi _{a, b} )}$ and $\ph ^{( \pi _{a, b, \overline{p}} )}$. 
Let 
\begin{align}
J(p) &= J_n (p; a, b, \overline{p} ) = %
\int_{0}^{1} {t^{a - 1} \{ 1 - p (1 - t) \} ^n \over \{ 1 - \overline{p} (1 - t) \} ^{n + a + b + 1}} dt \text{.} \label{eq:J} 
\end{align}
An upper bound on the risk difference $R(p, \ph ^{( \pi _{a, b, \overline{p}} )} ) - R(p, \ph ^{( \pi _{a, b} )} )$ is given in the following theorem, whose proof is at the end of this section. 

\begin{thm}
\label{thm:sufficient} 
The risk difference between $\ph ^{( \pi _{a, b} )}$ and $\ph ^{( \pi _{a, b, \overline{p}} )}$ satisfies 
\begin{align}
&{\{ R(p, \ph ^{( \pi _{a, b, \overline{p}} )} ) - R(p, \ph ^{( \pi _{a, b} )} ) \} / J(p) \over E_{p}^{X} [ 1 / I(X + a, n + a + b + 1, \overline{p} ) ]} \non \\
&\le (1 - p) \log \Big\{ 1 - {1 \over (1 - \overline{p} ) (n + a + b) J(p)} \Big\} + p \log \Big[ 1 + {1 \over \overline{p} (n + a + b)} \Big\{ 1 + {1 \over J(p)} \Big\} \Big] \text{,} \label{eq:thm_sufficient} 
\end{align}
where $J(p)$ is given by (\ref{eq:J}). 
\end{thm}

In order to use Theorem \ref{thm:sufficient}, we have to evaluate the integral $J(p)$ either numerically or analytically. 
First, we have 
\begin{align}
J(p) %
&= \int_{1}^{(1 - p) / (1 - \overline{p} )} \Big\{ {1 - p - (1 - \overline{p} ) u \over \overline{p} u - p} \Big\} ^{a - 1} {( \overline{p} u - p)^{a + b - 1} \over ( \overline{p} - p)^{a + b}} u^n du \non 
\end{align}
if $p \in (0, \overline{p} )$. 
Therefore, in the case of the uniform prior ($a = b = 1$), we can numerically evaluate $J(p)$ for all $p \in (0, \overline{p} ]$. 

Next,  $J(p) / (1 - p)$ times the right-hand side of (\ref{eq:thm_sufficient}) is bounded above by 
\begin{align}
J(0) \log \Big\{ 1 - {1 \over (1 - \overline{p} ) (n + a + b) J(0)} \Big\} + {p \over 1 - p} J(p) \log \Big[ 1 + {1 \over \overline{p} (n + a + b)} \Big\{ 1 + {1 \over J( \overline{p} )} \Big\} \Big] \text{.} \label{eq:sufficient_above} 
\end{align}
Since $\{ p / (1 - p) \} J(p) \le \{ \overline{p} / (1 - \overline{p} ) \} J(0)$, the risk difference is negative if 
\begin{align}
\log \Big\{ 1 - {1 \over (1 - \overline{p} ) (n + a + b) J(0)} \Big\} + {\overline{p} \over 1 - \overline{p}} \log \Big[ 1 + {1 \over \overline{p} (n + a + b)} \Big\{ 1 + {1 \over J( \overline{p} )} \Big\} \Big] < 0 \text{,} \non 
\end{align}
which is satisfied if $\overline{p}$ is sufficiently small. 
Also, this condition can be numerically checked since 
\begin{align}
J(0) &= %
I(a, n + a + b + 1, \overline{p} ) \quad \text{and} \quad J( \overline{p} ) = %
I(a, a + b + 1, \overline{p} ) \text{.} \non 
\end{align}

Finally, if $\overline{p} \le 1 / n$, then $\{ p / (1 - p) \} J(p)$ is nondecreasing in $p$ for $p \in (0, \overline{p} ]$ and, by (\ref{eq:sufficient_above}), a sufficient condition for $\ph ^{( \pi _{a, b, \overline{p}} )}$ to dominate $\ph ^{( \pi _{a, b} )}$ is that 
\begin{align}
- {1 \over (1 - \overline{p} ) (n + a + b)} + {\overline{p} \over 1 - \overline{p}} J( \overline{p} ) \log \Big[ 1 + {1 \over \overline{p} (n + a + b)} \Big\{ 1 + {1 \over J( \overline{p} )} \Big\} \Big] < 0 \text{.} \non 
\end{align}
The integral $J( \overline{p} )$ can be expressed in closed form as 
\begin{align}
J( \overline{p} ) &= I(a, a + 2, \overline{p} ) = {1 + I(a, a + 1, \overline{p} ) \over (1 - \overline{p} ) (a + 1)} = {1 + (1 - \overline{p} ) a \over (1 - \overline{p} )^2 a (a + 1)} \non 
\end{align}
by parts (iii), (ii), and (v) of Lemma \ref{lem:I} if $b = 1$ and as 
\begin{align}
J( \overline{p} ) &= I(1 / 2, 2, \overline{p} ) = {1 \over 1 - \overline{p}} \Big\{ 1 + {1 \over \sqrt{\overline{p} (1 - \overline{p} )}} \arctan {\sqrt{\overline{p} \over 1 - \overline{p}}} \Big\} \non 
\end{align}
by part (vi) of Lemma \ref{lem:I} if $a = b = 1 / 2$.

\begin{lem}
\label{lem:I} 
The function $I( \al , \ga , \overline{p} )$, $\ga > \al > 0$, has the following properties: 
\begin{itemize}
\item[{\rm{(i)}}]
$\al I( \al , \ga , \overline{p} ) = 1 + \overline{p} \ga I( \al + 1, \ga + 1, \overline{p} )$. 
\item[{\rm{(ii)}}]
$[1 + 1 / \{ ( \ga - \al ) I( \al , \ga , \overline{p} ) \} ] [1 + 1 / \{ \overline{p} \ga I( \al + 1, \ga + 1, \overline{p} ) \} ] = 1 + 1 / \{ \overline{p} ( \ga - \al ) I( \al + 1, \ga + 1, \overline{p} ) \} $. 
\item[{\rm{(iii)}}]
$1 + \overline{p} ( \ga - \al ) I( \al + 1, \ga + 1, \overline{p} ) = (1 - \overline{p} ) \al I( \al , \ga + 1, \overline{p} )$. 
\item[{\rm{(iv)}}]
$1 / I( \al + 1, \ga + 1, \overline{p} ) \le 1 + 1 / I( \al , \ga + 1, \overline{p} )$. 
\item[{\rm{(v)}}]
$I( \al , \al + 1, \overline{p} ) = 1 / \{ (1 - \overline{p} ) \al \} $. 
\item[{\rm{(vi)}}]
$I(1 / 2, 2, \overline{p} ) = \{ 1 / (1 - \overline{p} ) \} [1 + \{ 1 / \sqrt{\overline{p} (1 - \overline{p} )} \} \arctan \sqrt{\overline{p} / (1 - \overline{p} )}]$. 
\end{itemize}
\end{lem}

\noindent
{\bf Proof of Theorem \ref{thm:sufficient}.} \ \ Let $\De = R(p, \ph ^{( \pi _{a, b, \overline{p}} )} ) - R(p, \ph ^{( \pi _{a, b} )} )$ be the risk difference. 
Then, by (\ref{eq:ratio_p}) and (\ref{eq:ratio_q}), 
\begin{align}
\De &= E_{p}^{X} \Big[ - p \log \Big\{ 1 - {1 \over (X + a) I(X + a, n + a + b, \overline{p} )} \Big\} \non \\
&\quad - (1 - p) \log \Big\{ 1 + {1 \over (n - X + b) I(X + a, n + a + b, \overline{p} )} \Big\} \Big] \text{.} \label{eq:tsufficientp1} 
\end{align}
By part (i) of Lemma \ref{lem:I}, 
\begin{align}
- p \log \Big\{ 1 - {1 \over (X + a) I(X + a, n + a + b, \overline{p} )} \Big\} &= p \log \Big\{ 1 + {1 \over (X + a) I(X + a, n + a + b, \overline{p} ) - 1} \Big\} \non \\
&= p \log \Big\{ 1 + {1 \over \overline{p} (n + a + b) I(X + a + 1, n + a + b + 1, \overline{p} )} \Big\} \text{.} \label{eq:tsufficientp2} 
\end{align}
On the other hand, by parts (ii) and (iii) of Lemma \ref{lem:I}, 
\begin{align}
&- (1 - p) \log \Big\{ 1 + {1 \over (n - X + b) I(X + a, n + a + b, \overline{p} )} \Big\} \non \\
&= (1 - p) \log \Big[ \Big\{ 1 + {1 \over \overline{p} (n + a + b) I(X + a + 1, n + a + b + 1, \overline{p} )} \Big\} \non \\
&\quad / \Big\{ 1 + {1 \over \overline{p} (n - X + b) I(X + a + 1, n + a + b + 1, \overline{p} )} \Big\} \Big] \non \\
&= (1 - p) \log \Big\{ 1 - {X + a \over n + a + b} {1 \over 1 + \overline{p} (n - X + b) I(X + a + 1, n + a + b + 1, \overline{p} )} \Big\} \non \\
&= (1 - p) \log \Big\{ 1 - {1 \over (1 - \overline{p} ) (n + a + b) I(X + a, n + a + b + 1, \overline{p} )} \Big\} \text{.} \non 
\end{align}
Therefore, 
\begin{align}
\De &= p E_{p}^{X} \Big[ \log \Big\{ 1 + {1 \over \overline{p} (n + a + b) I(X + a + 1, n + a + b + 1, \overline{p} )} \Big\} \Big] \non \\
&\quad + (1 - p) E_{p}^{X} \Big[ \log \Big\{ 1 - {1 \over (1 - \overline{p} ) (n + a + b) I(X + a, n + a + b + 1, \overline{p} )} \Big\} \Big] \text{.} \non 
\end{align}

By part (iv) of Lemma \ref{lem:I}, 
\begin{align}
&E_{p}^{X} \Big[ \log \Big\{ 1 + {1 \over \overline{p} (n + a + b) I(X + a + 1, n + a + b + 1, \overline{p} )} \Big\} \non \\
&\le E_{p}^{X} \Big[ \log \Big\{ 1 + {1 \over \overline{p} (n + a + b)} + {1 \over \overline{p} (n + a + b) I(X + a, n + a + b + 1, \overline{p} )} \Big\} \Big] \text{.} \non 
\end{align}
Note that for any $M_1 , M_2 \ge 0$, $\ze \log (1 + M_1 + M_2 / \ze )$ is a nondecreasing concave function of $\ze $ for $\ze > 0$. 
Then, by the covariance inequality and Jensen's inequality, 
\begin{align}
&E_{p}^{X} \Big[ \log \Big\{ 1 + {1 \over \overline{p} (n + a + b)} + {1 \over \overline{p} (n + a + b) I(X + a, n + a + b + 1, \overline{p} )} \Big\} \Big] \non \\
&\le E_{p}^{X} \Big[ {1 \over I(X + a, n + a + b + 1, \overline{p} )} \Big] \non \\
&\quad \times E_{p}^{X} \Big[ I(X + a, n + a + b + 1, \overline{p} ) \log \Big\{ 1 + {1 \over \overline{p} (n + a + b)} + {1 \over \overline{p} (n + a + b) I(X + a, n + a + b + 1, \overline{p} )} \Big\} \Big] \non \\
&\le E_{p}^{X} \Big[ {1 \over I(X + a, n + a + b + 1, \overline{p} )} \Big] \non \\
&\quad \times E_{p}^{X} [ I(X + a, n + a + b + 1, \overline{p} ) ] \log \Big\{ 1 + {1 \over \overline{p} (n + a + b)} + {1 \over \overline{p} (n + a + b) E_{p}^{X} [ I(X + a, n + a + b + 1, \overline{p} ) ]} \Big\} \text{.} \non 
\end{align}
Similarly, since $\ze \log (1 - 1 / \ze )$ is a nondecreasing concave function of $\ze $ for $\ze > 1$, we have, by the covariance inequality and Jensen's inequality, 
\begin{align}
&E_{p}^{X} \Big[ \log \Big\{ 1 - {1 \over (1 - \overline{p} ) (n + a + b) I(X + a, n + a + b + 1, \overline{p} )} \Big\} \Big] \non \\
&\le E_{p}^{X} \Big[ {1 \over I(X + a, n + a + b + 1, \overline{p} )} \Big] \non \\
&\quad \times E_{p}^{X} \Big[ I(X + a, n + a + b + 1, \overline{p} ) \log \Big\{ 1 - {1 \over (1 - \overline{p} ) (n + a + b) I(X + a, n + a + b + 1, \overline{p} )} \Big\} \Big] \non \\
&\le E_{p}^{X} \Big[ {1 \over I(X + a, n + a + b + 1, \overline{p} )} \Big] \non \\
&\quad \times E_{p}^{X} [ I(X + a, n + a + b + 1, \overline{p} ) ] \log \Big\{ 1 - {1 \over (1 - \overline{p} ) (n + a + b) E_{p}^{X} [ I(X + a, n + a + b + 1, \overline{p} ) ]} \Big\} \text{.} \non 
\end{align}
Thus, since 
\begin{align}
E_{p}^{X} [ I(X + a, n + a + b + 1, \overline{p} ) ] &= \int_{0}^{1} {t^{a - 1} E_{p}^{X} [ t^{X} ] \over \{ 1 - \overline{p} (1 - t) \} ^{n + a + b + 1}} dt \non \\
&= \int_{0}^{1} {t^{a - 1} \{ 1 - p (1 - t) \} ^{n} \over \{ 1 - \overline{p} (1 - t) \} ^{n + a + b + 1}} dt = J(p) \text{,} \non 
\end{align}
we conclude that 
\begin{align}
&{\De / J(p) \over E_{p}^{X} [ 1 / I(X + a, n + a + b + 1, \overline{p} ) ]} \non \\
&\le p \log \Big\{ 1 + {1 \over \overline{p} (n + a + b)} + {1 \over \overline{p} (n + a + b) J(p)} \Big\} + (1 - p) \log \Big\{ 1 - {1 \over (1 - \overline{p} ) (n + a + b) J(p)} \Big\} \text{.} \non 
\end{align}
This completes the proof. 
\hfill$\Box$

\bigskip

The development in this section is analogous to that in Section 2 of Hamura and Kubokawa (2020), who used properties of the incomplete gamma function instead of the incomplete beta function in order to obtain a dominance condition in the Poisson case. 
However, the developments are not the same. 
First, both terms in (\ref{eq:loss_entropy}) are nonlinear, whereas only one term is nonlinear in the Poisson case. 
Second, as will be seen in Section \ref{sec:relation_Po}, the prior considered in their paper corresponds to $b = 1$. 
Therefore, an additional complication arises when we want to evaluate $J( \overline{p} )$ for the case of the Jeffreys prior ($a = b = 1 / 2$).

\subsection{Necessary conditions for domination}
\label{subsec:necessary}
It is shown in the following theorem that $\ph ^{( \pi _{a, b, \overline{p}} )}$ does not always dominate $\ph ^{( \pi _{a, b} )}$. 

\begin{thm}
\label{thm:necessary} 
If $\ph ^{( \pi _{a, b, \overline{p}} )}$ dominates $\ph ^{( \pi _{a, b} )}$, then $\overline{p} < (n + a) / (n + a + b)$. 
\end{thm}

\noindent
{\bf Proof%
.} \ \ By (\ref{eq:tsufficientp1}), (\ref{eq:tsufficientp2}), and part (ii) of Lemma \ref{lem:I}, 
\begin{align}
R( \overline{p} , \ph ^{( \pi _{a, b, \overline{p}} )} ) - R( \overline{p} , \ph ^{( \pi _{a, b} )} ) &= E_{\overline{p}}^{X} \Big[ \overline{p} \log \Big\{ 1 + {1 \over \overline{p} (n + a + b) I(X + a + 1, n + a + b + 1, \overline{p} )} \Big\} \non \\
&\quad - (1 - \overline{p} ) \log \Big\{ 1 + {1 \over (n - X + b) I(X + a, n + a + b, \overline{p} )} \Big\} \Big] \label{tnecessaryp1} \\
&= E_{\overline{p}}^{X} \Big[ \log \Big\{ 1 + {1 \over \overline{p} (n + a + b) I(X + a + 1, n + a + b + 1, \overline{p} )} \Big\} \non \\
&\quad - (1 - \overline{p} ) \log \Big\{ 1 + {1 \over \overline{p} (n - X + b) I(X + a + 1, n + a + b + 1, \overline{p} )} \Big\} \Big] \text{.} \non 
\end{align}
Therefore, since $(1 - \overline{p} ) \log (1 + \ze ) < \log \{ 1 + (1 - \overline{p} ) \ze \} $ for all $\ze > 0$, 
\begin{align}
R( \overline{p} , \ph ^{( \pi _{a, b, \overline{p}} )} ) - R( \overline{p} , \ph ^{( \pi _{a, b} )} ) &> E_{\overline{p}}^{X} \Big[ \log \Big\{ 1 + {1 \over \overline{p} (n + a + b) I(X + a + 1, n + a + b + 1, \overline{p} )} \Big\} \non \\
&\quad - \log \Big\{ 1 + {1 - \overline{p} \over \overline{p} (n - X + b) I(X + a + 1, n + a + b + 1, \overline{p} )} \Big\} \Big] \text{,} \non 
\end{align}
which is nonnegative if $1 / (n + a + b) \ge (1 - \overline{p} ) / b$ or $\overline{p} \ge (n + a) / (n + a + b)$. 
This proves Theorem \ref{thm:necessary}. 
\hfill$\Box$

\bigskip

Although the condition of Theorem \ref{thm:necessary} is not restrictive when $n$ is large, it is important that for any $a, b > 0$, the condition is not satisfied when $\overline{p}$ is sufficiently large. 
This is in contrast to the case of Theorem 2.2 of Hamura and Kubokawa (2020). 
The necessary condition of that theorem can be violated only when the hyperparameter $\be $ there satisfies $\log (1 + 1 / \be ) - 1 - 1 / (1 + \be ) > 0$. 

The following theorem gives a necessary condition which is similar to that of Theorem 2.2 of Hamura and Kubokawa (2020). 

\begin{thm}
\label{thm:necessary_2} 
Assume that $b = 1$. 
If $\ph ^{( \pi _{a, b, \overline{p}} )}$ dominates $\ph ^{( \pi _{a, b} )}$, then 
\begin{align}
\overline{p} \log \Big\{ 1 + {(1 - \overline{p} ) (a + 1) \over \overline{p} (n + a + 1)} \Big\} < (1 - \overline{p} ) \log {(n + a + 1) (1 - \overline{p} ^{n + 1} ) \over (n + 1) (1 - \overline{p} )} \text{.} \non 
\end{align}
\end{thm}

\noindent
{\bf Proof%
.} \ \ By (\ref{tnecessaryp1}), 
\begin{align}
R( \overline{p} , \ph ^{( \pi _{a, b, \overline{p}} )} ) - R( \overline{p} , \ph ^{( \pi _{a, b} )} ) &= \overline{p} E_{\overline{p}}^{X} \Big[ \log \Big\{ 1 + {1 \over \overline{p} (n + a + b) I(X + a + 1, n + a + b + 1, \overline{p} )} \Big\} \Big] \non \\
&\quad - (1 - \overline{p} ) E_{\overline{p}}^{X} \Big[ \log \Big\{ 1 + {1 \over (n - X + b) I(X + a, n + a + b, \overline{p} )} \Big\} \Big] \text{.} \non 
\end{align}
By Jensen's inequality, 
\begin{align}
&E_{\overline{p}}^{X} \Big[ \log \Big\{ 1 + {1 \over \overline{p} (n + a + b) I(X + a + 1, n + a + b + 1, \overline{p} )} \Big\} \Big] \non \\
&\ge \log \Big\{ 1 + {1 \over \overline{p} (n + a + b) E_{\overline{p}}^{X} [ I(X + a + 1, n + a + b + 1, \overline{p} ) ]} \Big\} \text{,} \non 
\end{align}
where 
\begin{align}
E_{\overline{p}}^{X} [ I(X + a + 1, n + a + b + 1, \overline{p} ) ] &= \int_{0}^{1} {t^{a} E_{\overline{p}}^{X} [ t^{X} ] \over \{ 1 - \overline{p} (1 - t) \} ^{n + a + b + 1}} dt \non \\
&= \int_{0}^{1} {t^{a} \{ 1 - \overline{p} (1 - t) \} ^{n} \over \{ 1 - \overline{p} (1 - t) \} ^{n + a + b + 1}} dt = I(a + 1, a + b + 1, \overline{p} ) \text{.} \non 
\end{align}
On the other hand, since 
\begin{align}
I(X + a, n + a + b, \overline{p} ) &= \int_{0}^{1} {t^{X + a - 1} \over \{ 1 - \overline{p} (1 - t) \} ^{n + a + b}} dt > \int_{0}^{1} t^{X + a - 1} dt = {1 \over X + a} \text{,} \non 
\end{align}
it follows that 
\begin{align}
&E_{\overline{p}}^{X} \Big[ \log \Big\{ 1 + {1 \over (n - X + b) I(X + a, n + a + b, \overline{p} )} \Big\} \Big] \non \\
&< E_{\overline{p}}^{X} \Big[ \log \Big( 1 + {X + a \over n - X + b} \Big) \Big] = E_{\overline{p}}^{X} \Big[ \log {n + a + b \over n - X + b} \Big] \le \log E_{\overline{p}}^{X} \Big[ {n + a + b \over n - X + b} \Big] \text{.} \non 
\end{align}
Therefore, 
\begin{align}
R( \overline{p} , \ph ^{( \pi _{a, b, \overline{p}} )} ) - R( \overline{p} , \ph ^{( \pi _{a, b} )} ) &> \overline{p} \log \Big\{ 1 + {1 \over \overline{p} (n + a + b) I(a + 1, a + b + 1, \overline{p} )} \Big\} - (1 - \overline{p} ) \log E_{\overline{p}}^{X} \Big[ {n + a + b \over n - X + b} \Big] \text{.} \non 
\end{align}

Now suppose that $b = 1$. 
Then, by part (v) of Lemma \ref{lem:I}, 
\begin{align}
I(a + 1, a + b + 1, \overline{p} ) = {1 \over (1 - \overline{p} ) (a + 1)} \text{.} \non 
\end{align}
By (3.4) of Chao and Strawderman (1972), 
\begin{align}
E_{\overline{p}}^{X} \Big[ {n + a + b \over n - X + b} \Big] &= {(n + a + 1) (1 - \overline{p} ^{n + 1} ) \over (n + 1) (1 - \overline{p} )} \text{.} \non 
\end{align}
Thus, 
\begin{align}
\overline{p} \log \Big\{ 1 + {(1 - \overline{p} ) (a + 1) \over \overline{p} (n + a + 1)} \Big\} < (1 - \overline{p} ) \log {(n + a + 1) (1 - \overline{p} ^{n + 1} ) \over (n + 1) (1 - \overline{p} )} \non 
\end{align}
if $\ph ^{( \pi _{a, b, \overline{p}} )}$ dominates $\ph ^{( \pi _{a, b} )}$. 
\hfill$\Box$

\section{The Case Where There are Both a Lower Bound Restriction and an Upper Bound Restriction}
\label{sec:both}
In this section, we assume $p \in [ \underline{p} , \overline{p} ]$ for known $0 < \underline{p} < \overline{p} < 1$ and compare $\pi _{a, b} (p) = p^{a - 1} (1 - p)^{b - 1}$ with $\pi _{a, b, \underline{p} , \overline{p}} (p) = \pi _{a, b} (p) 1_{[ \underline{p} , \overline{p} ]} (p)$. 
As in the previous section, the discussion focuses on the point estimation of $p$ on the basis of $X \sim {\rm{Bin}} (n, p)$ and we write $\ph ^{( \pi )} = \ph _{n}^{( \pi )} (X)$ for a prior $\pi (p)$ and $R(p, \de ) = %
E_{p}^{X} [ L(\de (X), p) ]$ for an estimator $\de (X)$ of $p$.

\subsection{Bayes estimators}
\label{subsec:both_estimators} 
The Bayes estimators of $p$ with respect to $\pi _{a, b} (p)$ and $\pi _{a, b, \underline{p} , \overline{p}} (p)$ are $\ph ^{( \pi _{a, b} )} = (X + a) / (n + a + b)$ and 
\begin{align}
\ph ^{( \pi _{a, b, \underline{p} , \overline{p}} )} = \frac{ \int_{\underline{p}}^{\overline{p}} p^{X + a} (1 - p)^{n - X + b - 1} dp }{ \int_{\underline{p}}^{\overline{p}} p^{X + a - 1} (1 - p)^{n - X + b - 1} dp } \text{,} \non 
\end{align}
respectively. 
As in Section \ref{subsec:estimators}, we have, by integration by parts, 
\begin{align}
\ph ^{( \pi _{a, b, \underline{p} , \overline{p}} )} &= \ph ^{( \pi _{a, b} )} - {A(X) \over n + a + b} \text{,} \non 
\end{align}
where 
\begin{align}
A(X) &= A_n (X; a, b, \underline{p} , \overline{p} ) = \Big[ p^{X + a} (1 - p)^{n - X + b} \Big] _{\underline{p}}^{\overline{p}} / \int_{\underline{p}}^{\overline{p}} p^{X + a - 1} (1 - p)^{n - X + b - 1} dp \text{.} \non 
\end{align}
It follows that 
\begin{align}
{\ph ^{( \pi _{a, b, \underline{p} , \overline{p}} )} \over \ph ^{( \pi _{a, b} )}} &= 1 - {A(X) \over X + a} \non 
\end{align}
and that 
\begin{align}
{1 - \ph ^{( \pi _{a, b, \underline{p} , \overline{p}} )} \over 1 - \ph ^{( \pi _{a, b} )}} &= 1 + {A(X) \over n - X + b} \text{.} \non 
\end{align}
However, in the present case, $A(X)$ is not always positive. 
If $a = b$ and $1 / 2 - \underline{p} = \overline{p} - 1 / 2 > 0$, then $A(X) \gtreqless 0$ if and only if $X \gtreqless n / 2$.

\subsection{The asymmetric case}
\label{subsec:asymmetric} 
Here, we compare the risk functions of $\ph ^{( \pi _{a, b} )}$ and $\ph ^{( \pi _{a, b, \underline{p} , \overline{p}} )}$. 

\begin{thm}
\label{thm:asymmetric} 
Suppose that 
\begin{align}
\overline{p} \le {a + 1 \over n + a + b + 1} \label{eq:assumption_asymmetric_0} 
\end{align}
and that 
\begin{align}
&{\overline{p} \over 1 - \overline{p}} \log {\underline{p} n + a + 1 \over \overline{p} (n + a + b)} + \log {(1 - \underline{p} ) n + b \over (1 - \overline{p} ) (n + a + b)} %
\le 0 \text{.} \label{eq:assumption_asymmetric} 
\end{align}
Then $\ph ^{( \pi _{a, b, \underline{p} , \overline{p}} )}$ dominates $\ph ^{( \pi _{a, b} )}$. 
\end{thm}

Clearly, both (\ref{eq:assumption_asymmetric_0}) and (\ref{eq:assumption_asymmetric}) are satisfied for sufficiently small $\overline{p}$. 
Also, the conditions are satisfied when $a$ is sufficiently large and $b$ is sufficiently small. 
In contrast, when $a = b$, the condition (\ref{eq:assumption_asymmetric_0}) implies that $\overline{p} \le 1 / 2$. 
Thus, Theorem \ref{thm:asymmetric} excludes the symmetric case of $a = b$ and $1 / 2 - \underline{p} = \overline{p} - 1 / 2 > 0$. 

A simpler condition for domination is given in the following corollary. 

\begin{cor}
\label{cor:asymmetric} 
Let $\overline{c} > \underline{c} > 0$ and assume $\underline{p} = \underline{c} / n$ and $\overline{p} = \overline{c} / n$. 
Suppose that $\overline{c} < a + 1$ and that $\overline{c} \log \{ ( \underline{c} + a + 1) / \overline{c} \} + \overline{c} - \underline{c} < a$. 
Then $\ph ^{( \pi _{a, b, \underline{p} , \overline{p}} )}$ dominates $\ph ^{( \pi _{a, b} )}$ for sufficiently large $n$. 
\end{cor}

\noindent
{\bf Proof%
.} \ \ Since $\overline{p} (n + a + b + 1) \to \overline{c} < a + 1$ as $n \to \infty $, the condition (\ref{eq:assumption_asymmetric_0}) is satisfied for large enough $n$. 
Since 
\begin{align}
&{1 \over 1 - \overline{p}} \log {\underline{p} n + a + 1 \over \overline{p} (n + a + b)} + {1 \over \overline{p}} \log {(1 - \underline{p} ) n + b \over (1 - \overline{p} ) (n + a + b)} \non \\
&= {1 \over 1 - \overline{c} / n} \log {\underline{c} + a + 1 \over \overline{c} \{ 1 + (a + b) / n \} } + {n \over \overline{c}} \log {n - \underline{c} + b \over n - \overline{c} + (1 - \overline{c} / n) (a + b)} \non \\
\log {\underline{c} + a + 1 \over \overline{c}} + {\overline{c} - \underline{c} - a \over \overline{c}} < 0 \non 
\end{align}
as $n \to \infty $, the condition (\ref{eq:assumption_asymmetric}) is satisfied for large enough $n$. 
\hfill$\Box$

\bigskip

We now prove Theorem \ref{thm:asymmetric}. 
Let $\underline{r} = \underline{p} / (1 - \underline{p} )$ and $\overline{r} = \overline{p} / (1 - \overline{p} )$ and let $\rho = \underline{r} / \overline{r} < 1$. 
Let, for $\ga > \al > 0$, 
\begin{align}
I( \al , \ga , \underline{p} , \overline{p} ) &= \int_{\rho }^{1} {t^{\al - 1} \over \{ 1 - \overline{p} (1 - t) \} ^{\ga }} dt \text{.} \non 
\end{align}
We use the following lemma. 

\begin{lem}
\label{lem:I_both} 
The function $I( \al , \ga , \underline{p} , \overline{p} )$, $\ga > \al > 0$, has the following properties: 
\begin{itemize}
\item[{\rm{(i)}}]
\begin{align}
{\al \over \ga } {\int_{\underline{p}}^{\overline{p}} p^{\al - 1} (1 - p)^{\ga - \al - 1} dp \over \int_{\underline{p}}^{\overline{p}} p^{\al } (1 - p)^{\ga - \al - 1} dp} = 1 + {1 \over \overline{p} \ga } \Big[ {t^{\al } \over \{ 1 - \overline{p} (1 - t) \} ^{\ga }} \Big] _{\rho }^{1} / I( \al + 1, \ga + 1, \underline{p} , \overline{p} ) \text{.} \non 
\end{align}
\item[{\rm{(ii)}}]
\begin{align}
{\ga - \al \over \ga } \frac{ \int_{\underline{p}}^{\overline{p}} p^{\al - 1} (1 - p)^{\ga - \al - 1} dp }{ \int_{\underline{p}}^{\overline{p}} p^{\al - 1} (1 - p)^{\ga - \al } dp } &= 1 - {1 \over (1 - \overline{p} ) \ga } \Big[ {t^{\al } \over \{ 1 - \overline{p} (1 - t) \} ^{\ga }} \Big] _{\rho }^{1} / I( \al , \ga + 1, \underline{p} , \overline{p} ) \text{.} \non 
\end{align}
\item[{\rm{(iii)}}]
\begin{align}
\Big[ {t^{\al } \over \{ 1 - \overline{p} (1 - t) \} ^{\ga }} \Big] _{\rho }^{1} / I( \al + 1, \ga + 1, \underline{p} , \overline{p} ) &< 1 + \Big[ {t^{\al } \over \{ 1 - \overline{p} (1 - t) \} ^{\ga }} \Big] _{\rho }^{1} / I( \al , \ga + 1, \underline{p} , \overline{p} ) \text{.} \non 
\end{align}
\end{itemize}
\end{lem}

\noindent
{\bf Proof of Theorem \ref{thm:asymmetric}.} \ \ By Lemma \ref{lem:I_both}, %
\begin{align}
&R(p, \ph ^{( \pi _{a, b, \underline{p} , \overline{p}} )} ) - R(p, \ph ^{( \pi _{a, b} )} ) \non \\
&= E_{p}^{X} \Big[ p \log \Big[ 1 + {1 \over \overline{p} (n + a + b)} \Big[ {t^{X + a} \over \{ 1 - \overline{p} (1 - t) \} ^{n + a + b}} \Big] _{\rho }^{1} / I(X + a + 1, n + a + b + 1, \underline{p} , \overline{p} ) \Big] \non \\
&\quad + (1 - p) \log \Big[ 1 - {1 \over (1 - \overline{p} ) (n + a + b)} \Big[ {t^{X + a} \over \{ 1 - \overline{p} (1 - t) \} ^{n + a + b}} \Big] _{\rho }^{1} / I(X + a, n + a + b + 1, \underline{p} , \overline{p} ) \Big] \Big] \non \\
&< %
E_{p}^{X} \Big[ h_n \Big( \Big[ {t^{X + a} \over \{ 1 - \overline{p} (1 - t) \} ^{n + a + b}} \Big] _{\rho }^{1} / I(X + a, n + a + b + 1, \underline{p} , \overline{p} ); p; a, b, \overline{p} \Big) \Big] \text{,} \non 
\end{align}
where 
\begin{align}
h_n ( \ze ; p; a, b, \overline{p} ) &= p \log \Big\{ 1 + {1 + \ze \over \overline{p} (n + a + b)} \Big\} + (1 - p) \log \Big\{ 1 - {\ze \over (1 - \overline{p} ) (n + a + b)} \Big\} \non 
\end{align}
for $- 1 - \overline{p} (n + a + b) < \ze < (1 - \overline{p} ) (n + a + b)$. 
For any $- (1 - \overline{p} ) \le \ze < (1 - \overline{p} ) (n + a + b)$, 
\begin{align}
{\pd h_n ( \ze ; p; a, b, \overline{p} ) \over \pd \ze } &= {p / \{ \overline{p} (n + a + b) \} \over 1 + (1 + \ze ) / \{ \overline{p} (n + a + b) \} } - {(1 - p) / \{ (1 - \overline{p} ) (n + a + b) \} \over 1 - \ze / \{ (1 - \overline{p} ) (n + a + b) \} } \non \\
&\le {1 / (n + a + b) \over 1 + (1 + \ze ) / \{ \overline{p} (n + a + b) \} } - {1 / (n + a + b) \over 1 - \ze / \{ (1 - \overline{p} ) (n + a + b) \} } \non \\
&\propto - \overline{p} \ze - (1 - \overline{p} ) (1 + \ze ) = - \ze - (1 - \overline{p} ) \le 0 \text{.} \non 
\end{align}
Note that 
\begin{align}
&\Big[ {t^{X + a} \over \{ 1 - \overline{p} (1 - t) \} ^{n + a + b}} \Big] _{\rho }^{1} / I(X + a, n + a + b + 1, \underline{p} , \overline{p} ) \non \\
&= \int_{\rho }^{1} \Big[ {(X + a) t^{X + a - 1} \over \{ 1 - \overline{p} (1 - t) \} ^{n + a + b}} - {\overline{p} (n + a + b) t^{X + a} \over \{ 1 - \overline{p} (1 - t) \} ^{n + a + b + 1}} \Big] dt / \int_{\rho }^{1} {t^{X + a - 1} \over \{ 1 - \overline{p} (1 - t) \} ^{n + a + b + 1}} dt \non \\
&\ge \inf_{t \in ( \rho , 1)} [(X + a) \{ 1 - \overline{p} (1 - t) \} - \overline{p} (n + a + b) t] = \inf_{t \in ( \rho , 1)} [(1 - \overline{p} ) (X + a) - \overline{p} (n - X + b) t] \non \\
&\ge (1 - \overline{p} ) (X + a) - \overline{p} (n - X + b) = X + a - \overline{p} (n + a + b) \ge a - \overline{p} (n + a + b) \text{.} \non 
\end{align}
Then, since $a - \overline{p} (n + a + b) \ge - (1 - \overline{p} )$ by (\ref{eq:assumption_asymmetric_0}), 
\begin{align}
R(p, \ph ^{( \pi _{a, b, \underline{p} , \overline{p}} )} ) - R(p, \ph ^{( \pi _{a, b} )} ) &< E_{p}^{X} [ h_n (X + a - \overline{p} (n + a + b); p; a, b, \overline{p} ) ] \non \\
&\le h_n (n p + a - \overline{p} (n + a + b); p; a, b, \overline{p} ) \non \\
&\le h_n (n \underline{p} + a - \overline{p} (n + a + b); p; a, b, \overline{p} ) \text{,} \non 
\end{align}
where the second inequality follows from Jensen's inequality. 
Furthermore, since $1 + n \underline{p} + a - \overline{p} (n + a + b) \ge 0$ by (\ref{eq:assumption_asymmetric_0}), 
\begin{align}
&h_n (n \underline{p} + a - \overline{p} (n + a + b); p; a, b, \overline{p} ) / (1 - p) \non \\
&= {p \over 1 - p} \log \Big\{ 1 + {1 + n \underline{p} + a - \overline{p} (n + a + b) \over \overline{p} (n + a + b)} \Big\} + \log \Big\{ 1 - {n \underline{p} + a - \overline{p} (n + a + b) \over (1 - \overline{p} ) (n + a + b)} \Big\} \non \\
&\le {\overline{p} \over 1 - \overline{p}} \log \Big\{ 1 + {1 + n \underline{p} + a - \overline{p} (n + a + b) \over \overline{p} (n + a + b)} \Big\} + \log \Big\{ 1 - {n \underline{p} + a - \overline{p} (n + a + b) \over (1 - \overline{p} ) (n + a + b)} \Big\} \text{,} \non 
\end{align}
which is nonpositive by (\ref{eq:assumption_asymmetric}). 
This completes the proof. 
\hfill$\Box$

\subsection{The symmetric case}
\label{subsec:symmetric} 
In this section, we assume that $a = b$ and that $1 / 2 - \underline{p} = \overline{p} - 1 / 2 > 0$. 
In this case, although there are different expressions for the risk difference, it is not clear whether they can be useful or not. 
Here, we take a direct approach. 

\begin{thm}
\label{thm:symmetric} 
Assume that $a = b$ and that $1 / 2 - \underline{p} = \overline{p} - 1 / 2 > 0$. 
Suppose that $a (a + 2)^3 > (a + 1)^4$. 
Then for $\underline{p}$ and $\overline{p}$ sufficiently close to $1 / 2$, the risk difference $R(p, \ph ^{( \pi _{a, b, \underline{p} , \overline{p}} )} ) - R(p, \ph ^{( \pi _{a, b} )} )$ is a convex function of $p$ %
and the maximum risk difference is negative. 
\end{thm}

The condition $a (a + 2)^3 > (a + 1)^4$ is satisfied when $a = 1$. 
This corresponds to the uniform prior. 
The condition is also satisfied when $a = 1 / 2$, which corresponds to the Jeffreys prior. 

Theorem \ref{thm:symmetric} assumes that $\underline{p}$ and $\overline{p}$ are close to $1 / 2$. 
In contrast, when $n = 1$, the maximum risk difference can be calculated exactly for any $( \underline{p} , \overline{p} )$. 

\begin{thm}
\label{thm:1n} 
Assume that $a = b$ and that $1 / 2 - \underline{p} = \overline{p} - 1 / 2 > 0$. 
Suppose that $n = 1$. 
Then the risk difference $R(p, \ph ^{( \pi _{a, b, \underline{p} , \overline{p}} )} ) - R(p, \ph ^{( \pi _{a, b} )} )$ is a convex function of $p$ and the maximum risk difference is 
\begin{align}
( \underline{p} ^2 + \overline{p} ^2 ) \log \Big\{ {1 + a \over 1 + 2 a} \frac{ \int_{\underline{p}}^{\overline{p}} p^{a} (1 - p)^{a - 1} dp }{ \int_{\underline{p}}^{\overline{p}} p^{1 + a} (1 - p)^{a - 1} dp } \Big\} + 2 \underline{p} \overline{p} \log \Big\{ {a \over 1 + 2 a} \frac{ \int_{\underline{p}}^{\overline{p}} p^{a - 1} (1 - p)^{a} dp }{ \int_{\underline{p}}^{\overline{p}} p^{a} (1 - p)^{a} dp } \Big\} \text{.} \label{eq:max_risk_1n} 
\end{align}
\end{thm}

If $a = 1$, (\ref{eq:max_risk_1n}) becomes 
\begin{align}
- ( \underline{p} ^2 + \overline{p} ^2 ) \log {\overline{p} ^3 - \underline{p} ^3 \over \overline{p} ^2 - \underline{p} ^2} - 2 \underline{p} \overline{p} \log \Big( 3 - 2 {\overline{p} ^3 - \underline{p} ^3 \over \overline{p} ^2 - \underline{p} ^2} \Big) \text{.} \non 
\end{align}
If $a = 1 / 2$, it becomes 
\begin{align}
&- ( \underline{p} ^2 + \overline{p} ^2 ) \log \Big\{ 1 - {2 \over 3} \frac{ \big[ u^{3 / 2} / (1 + u)^2 \big] _{\underline{r}}^{\overline{r}} }{ \big[ \arctan u \big] _{\sqrt{\underline{r}}}^{\sqrt{\overline{r}}} - \big[ u^{1 / 2} / (1 + u) \big] _{\underline{r}}^{\overline{r}} } \Big\} - 2 \underline{p} \overline{p} \log \Big\{ 1 + 2 \frac{ \big[ u^{3 / 2} / (1 + u)^2 \big] _{\underline{r}}^{\overline{r}} }{ \big[ \arctan u \big] _{\sqrt{\underline{r}}}^{\sqrt{\overline{r}}} - \big[ u^{1 / 2} / (1 + u) \big] _{\underline{r}}^{\overline{r}} } \Big\} \text{.} \non 
\end{align}

In the remainder of this section, we prove Theorems \ref{thm:symmetric} and \ref{thm:1n}. 
The following two lemmas are used in proving Theorem \ref{thm:symmetric}. 

\begin{lem}
\label{lem:partialp} 
Let $X \sim {\rm{Bin}} (n, p)$ and let $\varphi (X)$ be a function of $X$. 
Then 
\begin{align}
\Big( {\pd \over \pd p} \Big) ^2 \{ p E_{p}^{X} [ \varphi (X) ] \} &= {1 \over p} E_{p}^{X} [ X \{ (X + 1) \varphi (X) - 2 X \varphi (X - 1) + (X - 1) \varphi (X - 2) \} ] \text{.} \non 
\end{align}
\end{lem}

\begin{lem}
\label{lem:log_Jensen} 
Let $T$ be a random variable taking values in $(0, 1)$. 
Suppose that $T$ has mean $\mu \in (0, 1)$ and variance $\si ^2 \in (0, \infty )$. 
Then 
\begin{align}
E^T [ \log (1 - T) ] &\le \log (1 - \mu ) - \si ^2 / 2 \text{.} \non 
\end{align}
\end{lem}

\noindent
{\bf Proof of Theorem \ref{thm:symmetric}.} \ \ Since $a = b$ and $\underline{p} + \overline{p} = 1$, the risk difference between $\ph ^{( \pi _{a, b} )}$ and $\ph ^{( \pi _{a, b, \underline{p} , \overline{p}} )}$ can be written as 
\begin{align}
R(p, \ph ^{( \pi _{a, b, \underline{p} , \overline{p}} )} ) - R(p, \ph ^{( \pi _{a, b} )} ) &= E_{p}^{X} \Big[ p \log \Big\{ {X + a \over n + 2 a} \frac{ \int_{\underline{p}}^{\overline{p}} p^{X + a - 1} (1 - p)^{n - X + a - 1} dp }{ \int_{\underline{p}}^{\overline{p}} p^{X + a} (1 - p)^{n - X + a - 1} dp } \Big\} \non \\
&\quad + (1 - p) \log \Big\{ {n - X + a \over n + 2 a} \frac{ \int_{\underline{p}}^{\overline{p}} p^{X + a - 1} (1 - p)^{n - X + a - 1} dp }{ \int_{\underline{p}}^{\overline{p}} p^{X + a - 1} (1 - p)^{n - X + a} dp } \Big\} \Big] \non \\
&= p E_{p}^{X} [ h_n (X; a, \overline{p} ) ] + (1 - p) E_{1 - p}^{X} [ h_n (X; a, \overline{p} ) ] \text{,} \label{tsymmetricp1} 
\end{align}
where 
\begin{align}
h_n ( \ze ; a, \overline{p} ) &= \log \Big\{ {\ze + a \over n + 2 a} \frac{ \int_{\underline{p}}^{\overline{p}} p^{\ze + a - 1} (1 - p)^{n - \ze + a - 1} dp }{ \int_{\underline{p}}^{\overline{p}} p^{\ze + a} (1 - p)^{n - \ze + a - 1} dp } \Big\} \non 
\end{align}
for $- a < \ze < n + a$. 
By Lemma \ref{lem:partialp}, 
\begin{align}
&p \Big( {\pd \over \pd p} \Big) ^2 \{ p E_{p}^{X} [ h_n (X; a, \overline{p} ) ] \} \non \\
&= E_{p}^{X} [ X \{ (X + 1) h_n (X; a, \overline{p} ) - 2 X h_n (X - 1; a, \overline{p} ) + (X - 1) h_n (X - 2; a, \overline{p} ) \} ] \text{.} \label{tsymmetricp2} 
\end{align}

For any $- a < \ze < n + a$, 
\begin{align}
&{\pd h_n ( \ze ; a, \overline{p} ) \over \pd \ze } \non \\
&= {1 \over \ze + a} + \frac{ \int_{\underline{p}}^{\overline{p}} [ \log \{ p / (1 - p) \} ] p^{\ze + a - 1} (1 - p)^{n - \ze + a - 1} dp }{ \int_{\underline{p}}^{\overline{p}} p^{\ze + a - 1} (1 - p)^{n - \ze + a - 1} dp } - \frac{ \int_{\underline{p}}^{\overline{p}} [ \log \{ p / (1 - p) \} ] p^{\ze + a} (1 - p)^{n - \ze + a - 1} dp }{ \int_{\underline{p}}^{\overline{p}} p^{\ze + a} (1 - p)^{n - \ze + a - 1} dp } \non \\
&\ge {1 \over \ze + a} - 2 \log {\overline{p} \over \underline{p}} \non 
\end{align}
and 
\begin{align}
&{\pd ^2 h_n ( \ze ; a, \overline{p} ) \over {\pd \ze }^2} \non \\
&= - {1 \over ( \ze + a)^2} \non \\
&\quad + \frac{ \int_{\underline{p}}^{\overline{p}} [ \log \{ p / (1 - p) \} ]^2 p^{\ze + a - 1} (1 - p)^{n - \ze + a - 1} dp }{ \int_{\underline{p}}^{\overline{p}} p^{\ze + a - 1} (1 - p)^{n - \ze + a - 1} dp } - \Big( \frac{ \int_{\underline{p}}^{\overline{p}} [ \log \{ p / (1 - p) \} ] p^{\ze + a - 1} (1 - p)^{n - \ze + a - 1} dp }{ \int_{\underline{p}}^{\overline{p}} p^{\ze + a - 1} (1 - p)^{n - \ze + a - 1} dp } \Big) ^2 \non \\
&\quad - \frac{ \int_{\underline{p}}^{\overline{p}} [ \log \{ p / (1 - p) \} ]^2 p^{\ze + a} (1 - p)^{n - \ze + a - 1} dp }{ \int_{\underline{p}}^{\overline{p}} p^{\ze + a} (1 - p)^{n - \ze + a - 1} dp } + \Big( \frac{ \int_{\underline{p}}^{\overline{p}} [ \log \{ p / (1 - p) \} ] p^{\ze + a} (1 - p)^{n - \ze + a - 1} dp }{ \int_{\underline{p}}^{\overline{p}} p^{\ze + a} (1 - p)^{n - \ze + a - 1} dp } \Big) ^2 \non \\
&\ge - {1 \over ( \ze + a)^2} \non \\
&\quad - \frac{ \int_{\underline{p}}^{\overline{p}} [ \log \{ p / (1 - p) \} ]^2 p^{\ze + a} (1 - p)^{n - \ze + a - 1} dp }{ \int_{\underline{p}}^{\overline{p}} p^{\ze + a} (1 - p)^{n - \ze + a - 1} dp } + \Big( \frac{ \int_{\underline{p}}^{\overline{p}} [ \log \{ p / (1 - p) \} ] p^{\ze + a} (1 - p)^{n - \ze + a - 1} dp }{ \int_{\underline{p}}^{\overline{p}} p^{\ze + a} (1 - p)^{n - \ze + a - 1} dp } \Big) ^2 \non \\
&\ge - {1 \over ( \ze + a)^2} - \Big( \log {\overline{p} \over \underline{p}} \Big) ^2 \text{,} \non 
\end{align}
where the first inequality follows from the covariance inequality. 
In particular, for any $1 \le \ze \le n$, 
\begin{align}
{\pd ^2 \{ ( \ze + 1) h_n ( \ze ; a, \overline{p} ) \} \over {\pd \ze }^2} %
&= 2 {\pd h_n ( \ze ; a, \overline{p} ) \over \pd \ze } + ( \ze + 1) {\pd ^2 h_n ( \ze ; a, \overline{p} ) \over {\pd \ze }^2} \non \\
&\ge {2 \over \ze + a} - 4 \log {\overline{p} \over \underline{p}} - {\ze + 1 \over ( \ze + a)^2} - (n + 1) \Big( \log {\overline{p} \over \underline{p}} \Big) ^2 \non \\
&\ge {2 a \over (n + a)^2} - 4 \log {\overline{p} \over \underline{p}} - (n + 1) \Big( \log {\overline{p} \over \underline{p}} \Big) ^2 \text{,} \non 
\end{align}
which tends to $2 a / (n + a)^2 > 0$ as $\overline{p} \to 1 / 2$. 
This implies that for $\overline{p}$ sufficiently close to $1 / 2$, we have $(X + 1) h_n (X; a, \overline{p} ) - 2 X h_n (X - 1; a, \overline{p} ) + (X - 1) h_n (X - 2; a, \overline{p} ) \ge 0$ if $X \ge 3$. 
Next, for any $0 \le \ze \le 1$, 
\begin{align}
{\pd h_n ( \ze ; a, \overline{p} ) \over \pd \ze } &\ge {1 \over 1 + a} - 2 \log {\overline{p} \over \underline{p}} \to {1 \over 1 + a} > 0 \non 
\end{align}
as $\overline{p} \to 1 / 2$. 
Therefore, for the case that $X = 1$, we have $2 h_n (1; a, \overline{p} ) - 2 h_n (0; a, \overline{p} ) \ge 0$ for $\overline{p}$ sufficiently close to $1 / 2$. 
Finally, for the case of $X = 2$, 
\begin{align}
&3 h_n (2; a, \overline{p} ) - 4 h_n (1; a, \overline{p} ) + h_n (0; a, \overline{p} ) \ge \log {(2 + a)^3 a \over (1 + a)^4} - 4 \log {\overline{p} \over \underline{p}} \to \log {(2 + a)^3 a \over (1 + a)^4} > 0 \non 
\end{align}
as $\overline{p} \to 1 / 2$ by assumption. 
Thus, when $\overline{p}$ is sufficiently close to $1 / 2$, (\ref{tsymmetricp2}) is nonnegative and (\ref{tsymmetricp1}) is a convex function of $p$. 

Now, let $Z \sim {\rm{Bin}} (1, \overline{p} )$ and $W | Z \sim {\rm{Bin}} (n, Z \overline{p} + (1 - Z) (1 - \overline{p} ))$. 
Then for $\overline{p}$ sufficiently close to $1 / 2$, since $R(p, \ph ^{( \pi _{a, b, \underline{p} , \overline{p}} )} ) - R(p, \ph ^{( \pi _{a, b} )} )$ is convex in $p \in [ \underline{p} , \overline{p} ]$ and symmetric around $1 / 2$, 
\begin{align}
R(p, \ph ^{( \pi _{a, b, \underline{p} , \overline{p}} )} ) - R(p, \ph ^{( \pi _{a, b} )} ) &\le E_{\overline{p}}^{X} \Big[ \overline{p} \log \Big\{ {X + a \over n + 2 a} \frac{ \int_{\underline{p}}^{\overline{p}} p^{X + a - 1} (1 - p)^{n - X + a - 1} dp }{ \int_{\underline{p}}^{\overline{p}} p^{X + a} (1 - p)^{n - X + a - 1} dp } \Big\} \non \\
&\quad + (1 - \overline{p} ) \log \Big\{ {n - X + a \over n + 2 a} \frac{ \int_{\underline{p}}^{\overline{p}} p^{X + a - 1} (1 - p)^{n - X + a - 1} dp }{ \int_{\underline{p}}^{\overline{p}} p^{X + a - 1} (1 - p)^{n - X + a} dp } \Big\} \Big] \non \\
&= E_{\overline{p}}^{W} \Big[ \log \Big\{ {W + a \over n + 2 a} \frac{ \int_{\underline{p}}^{\overline{p}} p^{W + a - 1} (1 - p)^{n - W + a - 1} dp }{ \int_{\underline{p}}^{\overline{p}} p^{W + a} (1 - p)^{n - W + a - 1} dp } \Big\} \Big] \non \\
&\le E_{\overline{p}}^{W} \Big[ \log \Big( 1 - {n - W + a \over n + 2 a} \Big) \Big] + \log {1 \over \underline{p}} \text{.} \non 
\end{align}
Note that the marginal variance of $W$ is $n \overline{p} (1 - \overline{p} ) \{ 1 + n (2 \overline{p} - 1)^2 \} $. 
Then, by Lemma \ref{lem:log_Jensen}, 
\begin{align}
&E_{\overline{p}}^{W} \Big[ \log \Big( 1 - {n - W + a \over n + 2 a} \Big) \Big] + \log {1 \over \underline{p}} \non \\
&\le \log \Big[ 1 - {n \{ 1 - \overline{p} ^2 - (1 - \overline{p} )^2 \} + a \over n + 2 a} \Big] - {n \overline{p} (1 - \overline{p} ) \{ 1 + n (2 \overline{p} - 1)^2 \} \over 2 (n + 2 a)^2} + \log {1 \over \underline{p}} \text{,} \non 
\end{align}
which tends to $- n / \{ 8 (n + 2 a)^2 \} < 0$ as $\overline{p} \to 1 / 2$. 
This completes the proof. 
\hfill$\Box$

\bigskip

\noindent
{\bf Proof of Theorem \ref{thm:1n}.} \ \ By the proof of Theorem \ref{thm:symmetric}, we have 
\begin{align}
R(p, \ph ^{( \pi _{a, b, \underline{p} , \overline{p}} )} ) - R(p, \ph ^{( \pi _{a, b} )} ) &= p E_{p}^{X} [ h_1 (X; a, \overline{p} ) ] + (1 - p) E_{1 - p}^{X} [ h_1 (X; a, \overline{p} ) ] \non \\
&= \{ p^2 + (1 - p)^2 \} h_1 (1; a, \overline{p} ) + 2 p (1 - p) h_1 (0; a, \overline{p} ) \text{.} \non 
\end{align}
Note that 
\begin{align}
h_1 (1; a, \overline{p} ) - h_1 (0; a, \overline{p} ) &= \log \Big\{ {1 + a \over 1 + 2 a} \frac{ \int_{\underline{p}}^{\overline{p}} p^{a} (1 - p)^{a - 1} dp }{ \int_{\underline{p}}^{\overline{p}} p^{1 + a} (1 - p)^{a - 1} dp } \Big\} - \log \Big\{ {a \over 1 + 2 a} \frac{ \int_{\underline{p}}^{\overline{p}} p^{a - 1} (1 - p)^{a} dp }{ \int_{\underline{p}}^{\overline{p}} p^{a} (1 - p)^{a} dp } \Big\} \non \\
&= \log \frac{ (1 + a) \int_{\underline{p}}^{\overline{p}} p^{a} (1 - p)^{a} dp }{ a \int_{\underline{p}}^{\overline{p}} p^{1 + a} (1 - p)^{a - 1} dp } \non \\
&= \log \frac{ (1 + a) \int_{\underline{p}}^{\overline{p}} p^{a} (1 - p)^{a} dp }{ \big[ - p^{1 + a} (1 - p)^{a} \big] _{\underline{p}}^{\overline{p}} + (1 + a) \int_{\underline{p}}^{\overline{p}} p^{a} (1 - p)^{a} dp } \non \\
&= - \log \Big\{ 1 - \frac{ ( \underline{p} \overline{p} )^a ( \overline{p} - \underline{p} ) }{ (1 + a) \int_{\underline{p}}^{\overline{p}} p^{a} (1 - p)^{a} dp } \Big\} > 0 \text{.} \non 
\end{align}
Then 
\begin{align}
{\pd ^2 \{ R(p, \ph ^{( \pi _{a, b, \underline{p} , \overline{p}} )} ) - R(p, \ph ^{( \pi _{a, b} )} ) \} \over {\pd p}^2} &= 4 h_1 (1; a, \overline{p} ) - 4 h_1 (0; a, \overline{p} ) > 0 \non 
\end{align}
and the risk difference is a convex function of $p$. 
Therefore, 
\begin{align}
&\sup_{p \in [ \underline{p} , \overline{p} ]} \{ R(p, \ph ^{( \pi _{a, b, \underline{p} , \overline{p}} )} ) - R(p, \ph ^{( \pi _{a, b} )} ) \} = ( \underline{p} ^2 + \overline{p} ^2 ) h_1 (1; a, \overline{p} ) + 2 \underline{p} \overline{p} h_1 (0; a, \overline{p} ) \non \\
&= ( \underline{p} ^2 + \overline{p} ^2 ) \log \Big\{ {1 + a \over 1 + 2 a} \frac{ \int_{\underline{p}}^{\overline{p}} p^{a} (1 - p)^{a - 1} dp }{ \int_{\underline{p}}^{\overline{p}} p^{1 + a} (1 - p)^{a - 1} dp } \Big\} + 2 \underline{p} \overline{p} \log \Big\{ {a \over 1 + 2 a} \frac{ \int_{\underline{p}}^{\overline{p}} p^{a - 1} (1 - p)^{a} dp }{ \int_{\underline{p}}^{\overline{p}} p^{a} (1 - p)^{a} dp } \Big\} \non 
\end{align}
and this completes the proof. 
\hfill$\Box$

\section{Relation to Poisson Problems}
\label{sec:relation_Po}
We have considered Bayesian point estimation and predictive density estimation in the binomial case. 
In particular, for the cases where the probability parameter is %
restricted, %
we have compared the risk functions of the Bayes estimators based on the truncated and untruncated beta priors. 
We here note that our problems are related to Poisson problems as treated by Hamura and Kubokawa (2020). 
Although the derivation given below is rather informal, it will serve the purpose. 

Let $\Xt \sim {\rm{Po}} (r \la )$ and $\Yt \sim {\rm{Po}} (s \la )$ be independent Poisson variables for known $r, s > 0$ and unknown $\la \in (0, \infty )$ and let $\pit ( \la )$ be a prior for $\la $. 
Let $\pi _{K}^{( \pit )} (p) = \pit (K p)$ for $K > 0$. 
Then, as $n / r \sim \la / p \sim K \to \infty $, 
\begin{align}
n \ph _{n}^{( \pi _{K}^{( \pit )} )} (X) / r &= {n \over r} \frac{ \int_{0}^{1} p^{1 + X} (1 - p)^{n - X} \pit (K p) dp }{ \int_{0}^{1} p^{X} (1 - p)^{n - X} \pit (K p) dp } \non \\
&= \frac{ \int_{0}^{n / r} \la ^{1 + X} (1 - r \la / n)^{n - X} \pit (K r \la / n) d\la }{ \int_{0}^{n / r} \la ^{X} (1 - r \la / n)^{n - X} \pit (K r \la / n) d\la } \non \\
&\to \frac{ \int_{0}^{\infty } \la ^{1 + \Xt } e^{- r \la } \pit ( \la ) d\la }{ \int_{0}^{\infty } \la ^{\Xt } e^{- r \la } \pit ( \la ) d\la } = E_{\pit }^{\la | \Xt } [ \la | \Xt ] \text{,} \non 
\end{align}
which we denote by $\lah _{r}^{( \pit )} ( \Xt )$. 
Similarly, as $n / r \sim \la / p \sim l / s \sim K \to \infty $, 
\begin{align}
\fh _{l, n}^{( \pi _{K}^{( \pit )} )} (Y; X) &= \binom{l}{Y} \frac{ \int_{0}^{1} p^{Y + X} (1 - p)^{l - Y + n - X} \pit (K p) dp }{ \int_{0}^{1} p^{X} (1 - p)^{n - X} \pit (K p) dp } \non \\
&= \binom{l}{Y} \Big( {r \over n} \Big) ^Y \frac{ \int_{0}^{n / r} \la ^{Y + X} (1 - r \la / n)^{l - Y + n - X} \pit (K r \la / n) d\la }{ \int_{0}^{n / r} \la ^{X} (1 - r \la / n)^{n - X} \pit (K r \la / n) d\la } \non \\
&\to {s^{\Yt } \over \Yt !} \frac{ \int_{0}^{\infty } \la ^{\Yt + \Xt } e^{- (s + r) \la } \pit ( \la ) d\la }{ \int_{0}^{\infty } \la ^{\Xt } e^{- r \la } \pit ( \la ) d\la } = E_{\pit }^{\la | \Xt } \Big[ {(s \la )^{\Yt } \over \Yt !} e^{- s \la } \Big| \Xt \Big] \text{,} \non 
\end{align}
which is denoted by $\gh _{s, r}^{( \pit )} ( \Yt ; \Xt )$. 
Furthermore, 
\begin{align}
{n \over r} R_{n} (p, \ph _{n}^{( \pi _{K}^{( \pit )} )} ) &= E_{p}^{X} \Big[ (1 - p) \log {(1 - p)^{n / r} \over \{ 1 - \ph _{n}^{( \pi _{K}^{( \pit )} )} (X) \} ^{n / r}} + {n p \over r} \log {n p  / r \over n \ph _{n}^{( \pi _{K}^{( \pit )} )} (X) / r} \Big] \non \\
&\to %
E_{\la }^{\Xt } \Big[ \lah _{r}^{( \pit )} ( \Xt ) - \la - \la \log {\lah _{r}^{( \pit )} ( \Xt ) \over \la } \Big] \text{,} \non 
\end{align}
which is denoted by $\widetilde{R} _r ( \la , \lah _{r}^{( \pit )} )$, and 
\begin{align}
R_{l, n} (p, \fh _{l, n}^{( \pi _{K}^{( \pit )} )} ) &= E_{p}^{(Y, X)} \Big[ \log {f(Y | l, p) \over \fh _{l, n}^{( \pi _{K}^{( \pit )} )} (Y ; X)} \Big] \non \\
&\to E_{\la }^{( \Yt , \Xt )} \Big[ \log {g( \Yt | s, \la ) \over \gh _{s, r}^{( \pit )} ( \Yt ; \Xt )} \Big] \text{,} \non 
\end{align}
where $g( \Yt | s, \la ) = \{ (s \la )^{\Yt } / ( \Yt !) \} e^{- s \la }$, and also 
\begin{align}
R_{l, n} (p, \fh _{l, n}^{( \pi _{K}^{( \pit )} )} ) &= \sum_{i = 0}^{l - 1} R_{n + i} (p, \ph _{n + i}^{( \pi _{K}^{( \pit )} )} ) \approx \sum_{i = 0}^{l - 1} {s \over l} \widetilde{R} _{r + i s / l} ( \la , \lah _{r + i s / l}^{( \pit )} ) \approx \int_{r}^{s + r} \widetilde{R} _{\ta } ( \la , \lah _{\ta }^{( \pit )} ) d\ta \text{.} \non 
\end{align}
Thus, the Bayesian procedures with respect to $\pi _{K}^{( \pit )} (p)$ and $X$ are asymptotically equivalent to those with respect to $\pit ( \la )$ and $\Xt $. 
Finally, for $\overline{\la } > 0$, the restriction $\la \in (0, \overline{\la } ]$ corresponds to $p \in (0, \overline{p} ] = (0, r \overline{\la } / n]$ (namely $\overline{p} \sim \overline{\la } / K$) and the priors $\pit _a ( \la ) = \la ^{a - 1}$ and $\pit _{a, \overline{\la }} ( \la ) = \pit _a ( \la ) 1_{(0, \overline{\la } ]} ( \la )$ correspond to $\pi _{K}^{( \pit _a )} (p) \propto \pi _{a, 1} (p)$ and $\pi _{K}^{( \pit _{a, \overline{\la }} )} (p) \propto \pi _{a, 1, \overline{\la } / K} (p)$, respectively. 
The restriction considered in Corollary \ref{cor:asymmetric} corresponds to $\la \in [ \underline{c} / r, \overline{c} / r]$ and this case was not treated by Hamura and Kubokawa (2020).

\section{Numerical Studies}
\label{sec:sim}
We first investigate %
the numerical performance of the risk functions of $\ph ^{( \pi _{a, b} )}$ and $\ph ^{( \pi _{a, b, \overline{p}} )}$ considered in Section \ref{sec:domination}. 
We set $a = b = 1$ since then the integral $J(p)$ appearing in (\ref{eq:thm_sufficient}) can be expressed in closed form. 
We consider the cases $n = 1$, $n = 5$, and $n = 9$. 
For each of these cases, we consider the cases $\overline{p} = 0.1$, $\overline{p} = 0.2$, $\overline{p} = 0.3$, and $\overline{p} = 0.4$. 
We note that since the sample space is finite (consisting of $n + 1$ points) in the binomial case, we do not need to use Monte Carlo simulation to approximate the risk functions. 

\begin{figure}
\centering
\includegraphics[width = 16cm]{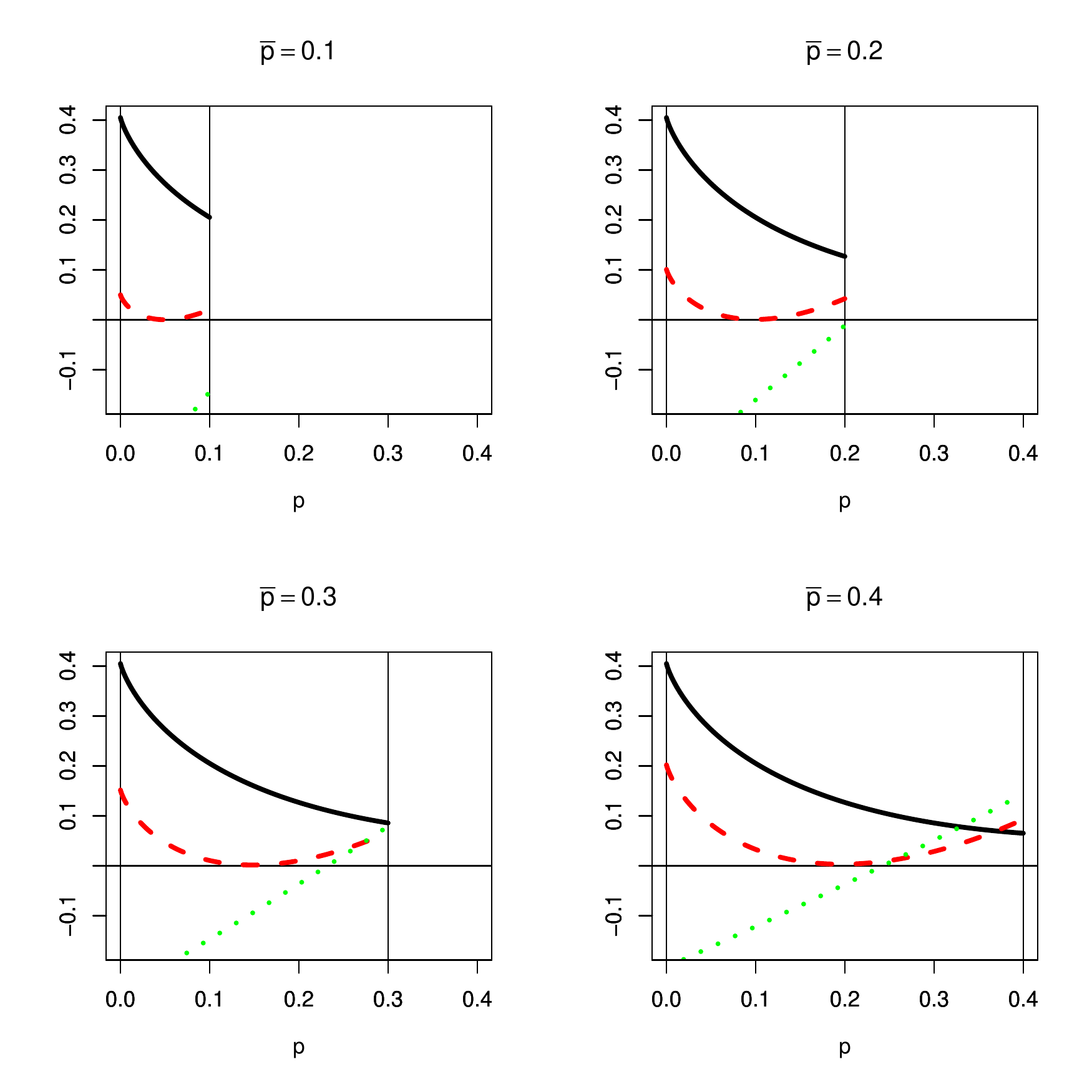}
\caption{Risks of the estimators $\ph ^{( \pi _{a, b} )}$ and $\ph ^{( \pi _{a, b, \overline{p}} )}$ with $a = b = 1$ for the cases $\overline{p} = 0.1$, $\overline{p} = 0.2$, $\overline{p} = 0.3$, and $\overline{p} = 0.4$. 
We set $n = 1$. 
The vertical lines show the boundary of $(0, \overline{p} ]$. 
The solid black and dashed red curves correspond to $\ph ^{( \pi _{a, b} )}$ and $\ph ^{( \pi _{a, b, \overline{p}} )}$, respectively. 
The dotted green curves show the right-hand side of (\ref{eq:thm_sufficient}). }
\label{fig:1n}
\end{figure}

\begin{figure}
\centering
\includegraphics[width = 16cm]{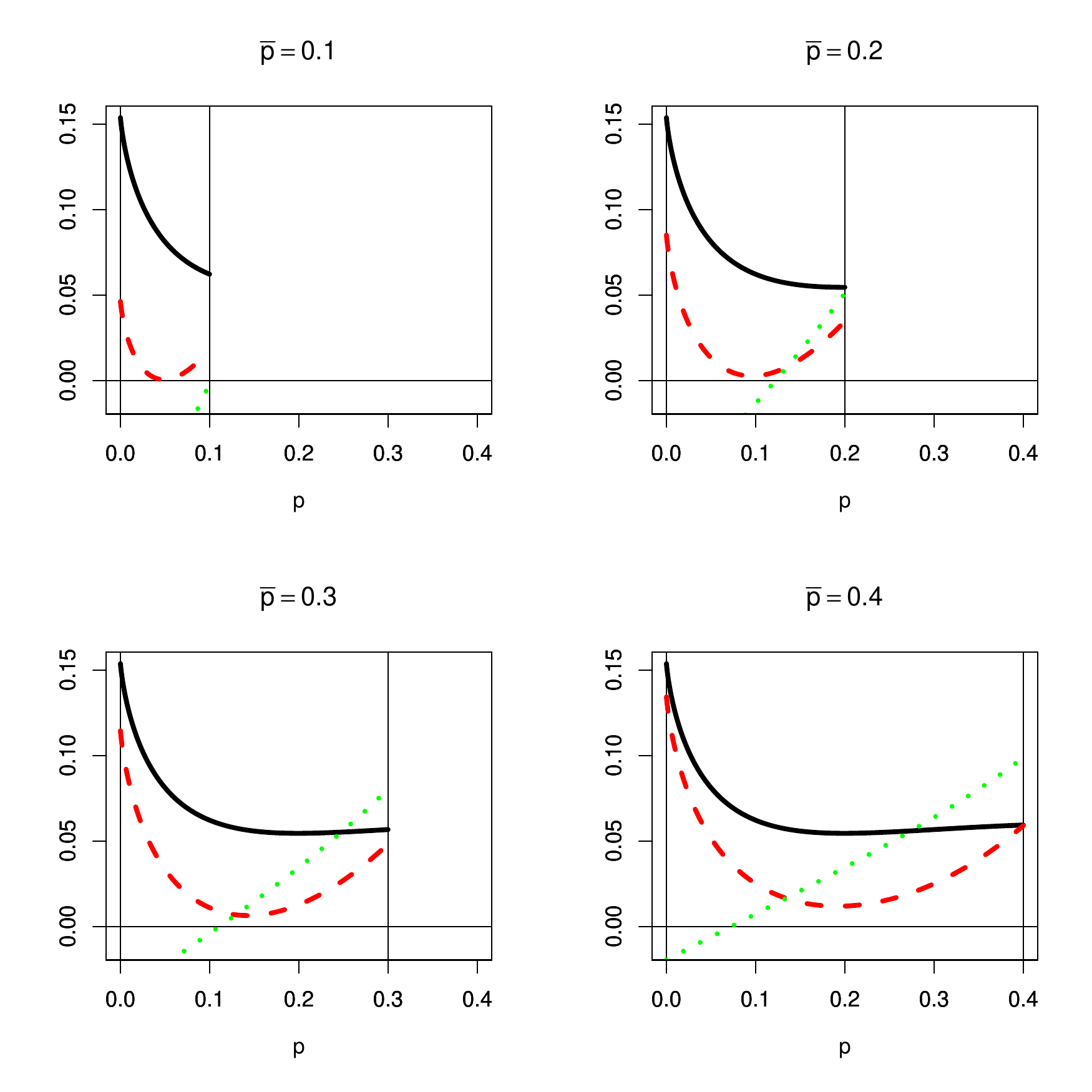}
\caption{Risks of the estimators $\ph ^{( \pi _{a, b} )}$ and $\ph ^{( \pi _{a, b, \overline{p}} )}$ with $a = b = 1$ for the cases $\overline{p} = 0.1$, $\overline{p} = 0.2$, $\overline{p} = 0.3$, and $\overline{p} = 0.4$. 
We set $n = 5$. 
The vertical lines show the boundary of $(0, \overline{p} ]$. 
The solid black and dashed red curves correspond to $\ph ^{( \pi _{a, b} )}$ and $\ph ^{( \pi _{a, b, \overline{p}} )}$, respectively. 
The dotted green curves show the right-hand side of (\ref{eq:thm_sufficient}). }
\label{fig:5n}
\end{figure}

\begin{figure}
\centering
\includegraphics[width = 16cm]{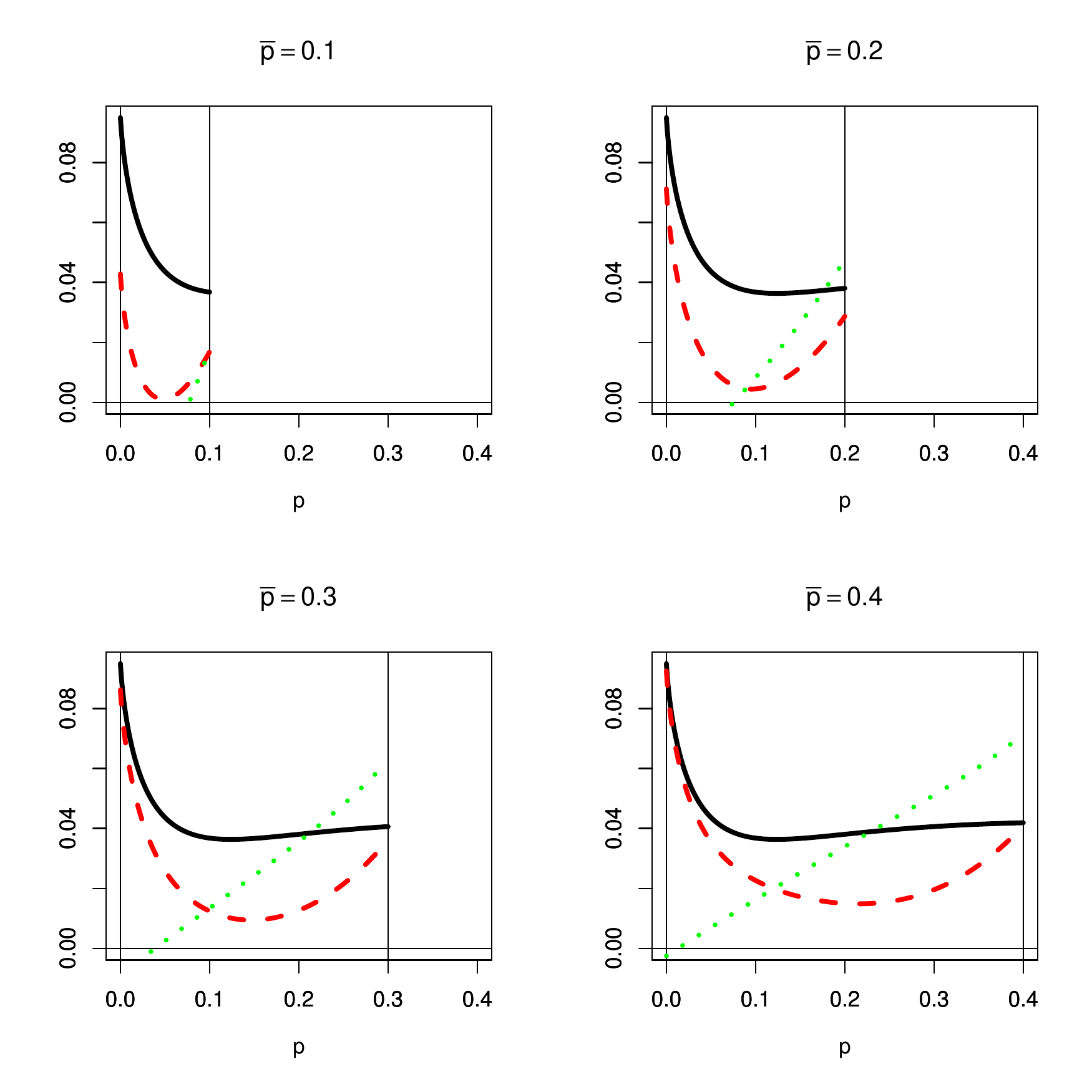}
\caption{Risks of the estimators $\ph ^{( \pi _{a, b} )}$ and $\ph ^{( \pi _{a, b, \overline{p}} )}$ with $a = b = 1$ for the cases $\overline{p} = 0.1$, $\overline{p} = 0.2$, $\overline{p} = 0.3$, and $\overline{p} = 0.4$. 
We set $n = 9$. 
The vertical lines show the boundary of $(0, \overline{p} ]$. 
The solid black and dashed red curves correspond to $\ph ^{( \pi _{a, b} )}$ and $\ph ^{( \pi _{a, b, \overline{p}} )}$, respectively. 
The dotted green curves show the right-hand side of (\ref{eq:thm_sufficient}). }
\label{fig:9n}
\end{figure}

The results are illustrated in Figures \ref{fig:1n}, \ref{fig:5n}, and \ref{fig:9n}. 
Overall, $\ph ^{( \pi _{a, b, \overline{p}} )}$ performs better than $\ph ^{( \pi _{a, b} )}$. 
The right-hand side of (\ref{eq:thm_sufficient}), which is an upper bound on the standardized risk difference, is negative when $n$ and $\overline{p}$ are small. 
However, when $n$ and $\overline{p}$ are large, the right-hand side of (\ref{eq:thm_sufficient}) tends to be large and therefore the sufficient conditions for domination will be restrictive. 

\begin{figure}
\centering
\includegraphics[width = 16cm]{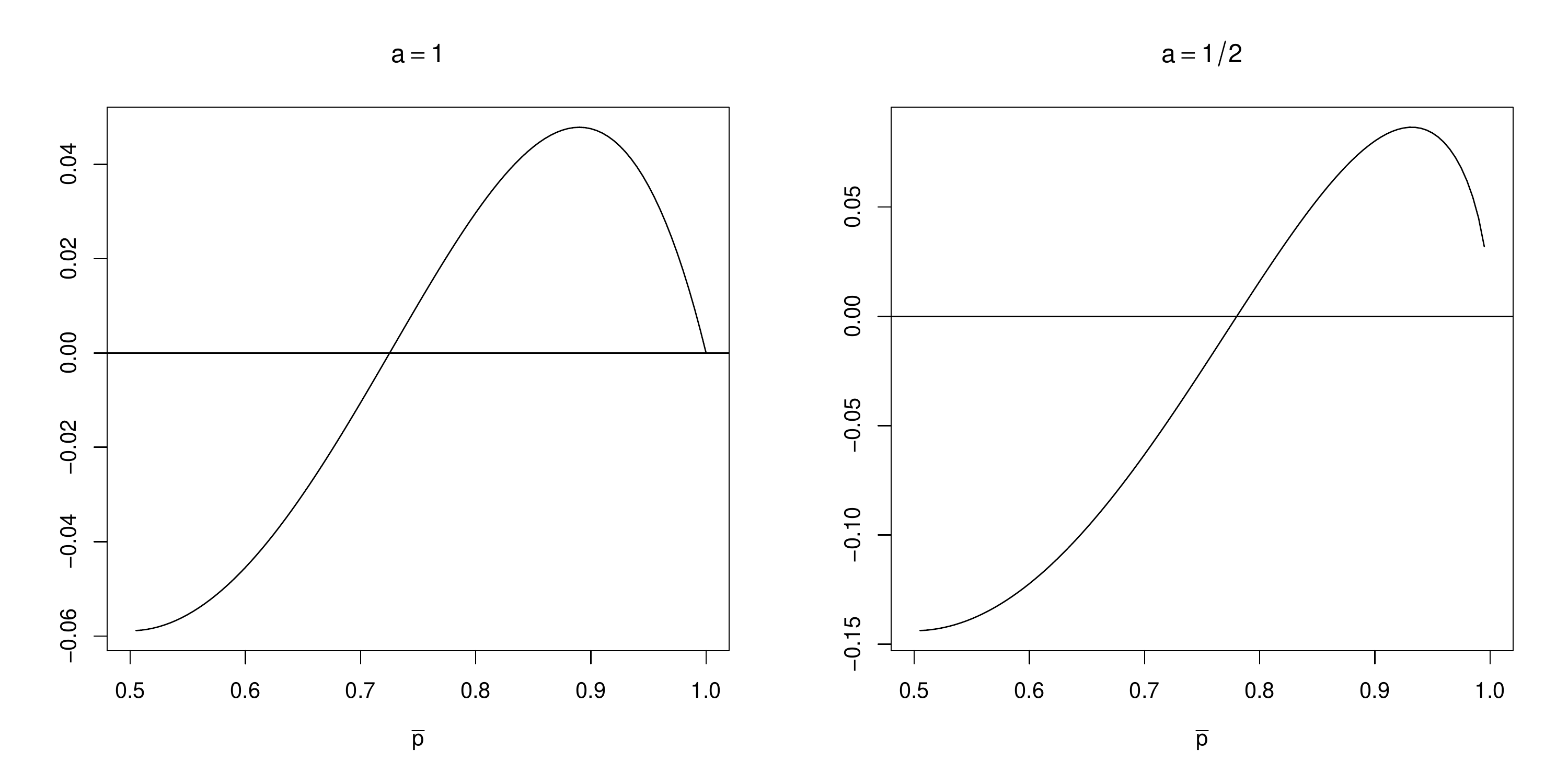}
\caption{The maximum risk difference (\ref{eq:max_risk_1n}) for the cases $a = 1$ and $a = 1 / 2$. }
\label{fig:maxRD_1n}
\end{figure}

Next, Figure \ref{fig:maxRD_1n} shows the maximum risk difference (\ref{eq:max_risk_1n}) for $1 / 2 < \overline{p} = 1 - \underline{p} < 1$ for each of the cases $a = b = 1$ and $a = b = 1 / 2$. 
If $n = 1$, $a = b = 1$, and $1 / 2 - \underline{p} = \overline{p} - 1 / 2 > 0$, then $\ph ^{( \pi _{a, b, \underline{p} , \overline{p}} )}$ dominates $\ph ^{( \pi _{a, b} )}$ if and only if $\overline{p} \lessapprox 0.725$. 
If $n = 1$, $a = b = 1 / 2$, and $1 / 2 - \underline{p} = \overline{p} - 1 / 2 > 0$, then $\ph ^{( \pi _{a, b, \underline{p} , \overline{p}} )}$ dominates $\ph ^{( \pi _{a, b} )}$ if and only if $\overline{p} \lessapprox 0.775$.

\section{Appendix}
Here we prove Lemmas \ref{lem:I}, \ref{lem:I_both}, \ref{lem:partialp}, and \ref{lem:log_Jensen}. 

\bigskip

\noindent
{\bf Proof of Lemma \ref{lem:I}.} \ \ For part (i), by integration by parts, we have 
\begin{align}
\al I( \al , \ga , \overline{p} ) &= \al \int_{0}^{1} {t^{\al - 1} \over \{ 1 - \overline{p} (1 - t) \} ^{\ga }} dt \non \\
&= \Big[ {t^{\al } \over \{ 1 - \overline{p} (1 - t) \} ^{\ga }} \Big] _{0}^{1} + \overline{p} \ga \int_{0}^{1} {t^{\al } \over \{ 1 - \overline{p} (1 - t) \} ^{\ga + 1}} dt = 1 + \overline{p} \ga I( \al + 1, \ga + 1, \overline{p} ) \text{.} \non 
\end{align}
Part (ii) follows from part (i): 
\begin{align}
&\Big\{ 1 + {1 \over ( \ga - \al ) I( \al , \ga , \overline{p} )} \Big\} \Big\{ 1 + {1 \over \overline{p} \ga I( \al + 1, \ga + 1, \overline{p} )} \Big\} \non \\
&= \Big\{ 1 + {\al / ( \ga - \al ) \over 1 + \overline{p} \ga I( \al + 1, \ga + 1, \overline{p} )} \Big\} {1 + \overline{p} \ga I( \al + 1, \ga + 1, \overline{p} ) \over \overline{p} \ga I( \al + 1, \ga + 1, \overline{p} )} \non \\
&= \Big\{ 1 + \overline{p} \ga I( \al + 1, \ga + 1, \overline{p} ) + {\al \over \ga - \al } \Big\} {1 \over \overline{p} \ga I( \al + 1, \ga + 1, \overline{p} )} \non \\
&= 1 + {1 \over \overline{p} ( \ga - \al ) I( \al + 1, \ga + 1, \overline{p} )} \text{.} \non 
\end{align}
For part (iii), note that by definition 
\begin{align}
I( \al , \ga , \overline{p} ) &= (1 - \overline{p} ) I( \al , \ga + 1, \overline{p} ) + \overline{p} I( \al + 1, \ga + 1, \overline{p} ) \text{.} \non 
\end{align}
Then, by part (i), 
\begin{align}
1 + \overline{p} \ga I( \al + 1, \ga + 1, \overline{p} ) &= (1 - \overline{p} ) \al I( \al , \ga + 1, \overline{p} ) + \overline{p} \al I( \al + 1, \ga + 1, \overline{p} ) \non 
\end{align}
or 
\begin{align}
1 + \overline{p} ( \ga - \al ) I( \al + 1, \ga + 1, \overline{p} ) = (1 - \overline{p} ) \al I( \al , \ga + 1, \overline{p} ) \text{.} \non 
\end{align}
For part (iv), let, for $t \in (0, 1)$, $h(t; \al , \ga , \overline{p} ) = t^{\al - 1} / \{ 1 - \overline{p} (1 - t) \} ^{\ga + 1} $. 
Then, by the covariance inequality, 
\begin{align}
{I( \al + 1, \ga + 2, \overline{p} ) \over I( \al , \ga + 1, \overline{p} )} &= \int_{0}^{1} {t \over 1 - \overline{p} + \overline{p} t} h(t; \al , \ga , \overline{p} ) dt / \int_{0}^{1} h(t; \al , \ga , \overline{p} ) dt \non \\
&\le \int_{0}^{1} {t^2 \over 1 - \overline{p} + \overline{p} t} h(t; \al , \ga , \overline{p} ) dt / \int_{0}^{1} t h(t; \al , \ga , \overline{p} ) dt = {I( \al + 2, \ga + 2, \overline{p} ) \over I( \al + 1, \ga + 1, \overline{p} )} \text{.} \non 
\end{align}
Therefore, by part (i), 
\begin{align}
{1 \over \overline{p} ( \ga + 1)} \Big\{ \al - {1 \over I( \al , \ga + 1, \overline{p} )} \Big\} %
{1 \over \overline{p} ( \ga + 1)} \Big\{ \al + 1 - {1 \over I( \al + 1, \ga + 1, \overline{p} )} \Big\} \non 
\end{align}
or 
\begin{align}
{1 \over I( \al + 1, \ga + 1, \overline{p} )} \le 1 + {1 \over I( \al , \ga + 1, \overline{p} )} \text{.} \non 
\end{align}
For part (v), 
\begin{align}
I( \al , \al + 1, \overline{p} ) &= \int_{0}^{1} {t^{\al - 1} \over \{ 1 - \overline{p} (1 - t) \} ^{\al + 1}} dt \non \\
&= \int_{0}^{\infty } {1 \over (1 + u)^{\al - 1}} {1 \over \{ 1 - \overline{p} u / (1 + u) \} ^{a + 1}} {1 \over (1 + u)^2} du = %
{1 \over (1 - \overline{p} ) \al } \text{.} \non 
\end{align}
For part (vi), 
\begin{align}
I(1 / 2, 2, \overline{p} ) &= \int_{0}^{1} {t^{1 / 2 - 1} \over \{ 1 - \overline{p} (1 - t) \} ^{2}} dt \non \\
&= \int_{0}^{1} {1 \over \tilde{t}} {1 \over \{ 1 - \overline{p} (1 - \tilde{t} ^2 ) \} ^{2}} 2 \tilde{t} d\tilde{t} = {2 \over \overline{p} ^2} \int_{0}^{1} {1 \over \{ (1 - \overline{p} ) / \overline{p} + \tilde{t} ^2 \} ^{2}} d\tilde{t} \non \\
&= {1 \over 1 - \overline{p}} \Big\{ 1 + {1 \over \sqrt{\overline{p} (1 - \overline{p} )}} \arctan {\sqrt{\overline{p} \over 1 - \overline{p}}} \Big\} \text{.} \non 
\end{align}
This completes the proof. 
\hfill$\Box$

\bigskip

\noindent
{\bf Proof of Lemma \ref{lem:I_both}.} \ \ For part (i), 
\begin{align}
\int_{\underline{p}}^{\overline{p}} p^{\al - 1} (1 - p)^{\ga - \al - 1} dp = \int_{\underline{r}}^{\overline{r}} {u^{\al - 1} \over (1 + u)^{\ga }} du = \overline{r} ^{\al } \int_{\rho }^{1} {t^{\al - 1} \over (1 + \overline{r} t)^{\ga }} dt \text{.} \non 
\end{align}
By integration by parts, 
\begin{align}
\al \int_{\rho }^{1} {t^{\al - 1} \over (1 + \overline{r} t)^{\ga }} dt &= \Big[ {t^{\al } \over (1 + \overline{r} t)^{\ga }} \Big] _{\rho }^{1} + \overline{r} \ga \int_{\rho }^{1} {t^{\al } \over (1 + \overline{r} t)^{\ga + 1}} dt \text{.} \label{lI_bothp1} 
\end{align}
Therefore, 
\begin{align}
\frac{ \int_{\underline{p}}^{\overline{p}} p^{\al - 1} (1 - p)^{\ga - \al - 1} dp }{ \int_{\underline{p}}^{\overline{p}} p^{\al } (1 - p)^{\ga - \al - 1} dp } &= {1 \over \overline{r} \al } \Big\{ \Big[ {t^{\al } \over (1 + \overline{r} t)^{\ga }} \Big] _{\rho }^{1} / \int_{\rho }^{1} {t^{\al } \over (1 + \overline{r} t)^{\ga + 1}} dt + \overline{r} \ga \Big\} \non \\
&= {\ga \over \al } \Big[ 1 + {1 \over \overline{r} \ga } {1 \over 1 - \overline{p}} \Big[ {t^{\al } \over \{ 1 - \overline{p} (1 - t) \} ^{\ga }} \Big] _{\rho }^{1} / \int_{\rho }^{1} {t^{\al } \over \{ 1 - \overline{p} (1 - t) \} ^{\ga + 1}} dt \Big] \non \\
&= {\ga \over \al } \Big[ 1 + {1 \over \overline{p} \ga } \Big[ {t^{\al } \over \{ 1 - \overline{p} (1 - t) \} ^{\ga }} \Big] _{\rho }^{1} / I( \al + 1, \ga + 1, \underline{p} , \overline{p} ) \Big] \text{.} \non 
\end{align}
For part (ii), it follows from part (i) that 
\begin{align}
\frac{ \int_{\underline{p}}^{\overline{p}} p^{\al - 1} (1 - p)^{\ga - \al - 1} dp }{ \int_{\underline{p}}^{\overline{p}} p^{\al - 1} (1 - p)^{\ga - \al } dp } &= \frac{ \int_{1 - \overline{p}}^{1 - \underline{p}} q^{\ga - \al - 1} (1 - q)^{\al - 1} dq }{ \int_{1 - \overline{p}}^{1 - \underline{p}} q^{\ga - \al } (1 - q)^{\al - 1} dq } \non \\
&= {\ga \over \ga - \al } \Big[ 1 + {1 \over (1 - \underline{p} ) \ga } \Big[ {t^{\ga - \al } \over \{ 1 - (1 - \underline{p} ) (1 - t) \} ^{\ga }} \Big] _{\rho }^{1} / \int_{\rho }^{1} {t^{\ga - \al } \over \{ 1 - (1 - \underline{p} ) (1 - t) \} ^{\ga + 1}} dt \Big] \text{.} \non 
\end{align}
Note that 
\begin{align}
&\Big[ {t^{\ga - \al } \over \{ 1 - (1 - \underline{p} ) (1 - t) \} ^{\ga }} \Big] _{\rho }^{1} / \int_{\rho }^{1} {t^{\ga - \al } \over \{ 1 - (1 - \underline{p} ) (1 - t) \} ^{\ga + 1}} dt \non \\
&= \Big[ {(1 / v)^{\ga - \al } \over \{ 1 - (1 - \underline{p} ) (1 - 1 / v) \} ^{\ga }} \Big] _{1 / \rho }^{1} / \int_{1}^{1 / \rho } {1 \over v^{\ga - \al }} {1 \over \{ 1 - (1 - \underline{p} ) (1 - 1 / v) \} ^{\ga + 1}} {1 \over v^2} dv \non \\
&= \Big[ {v^{\al } \over (1 - \underline{p} + \underline{p} v)^{\ga }} \Big] _{1 / \rho }^{1} / \int_{1}^{1 / \rho } {v^{\al - 1} \over (1 - \underline{p} + \underline{p} v)^{\ga + 1}} dv \non \\
&= \Big[ {(t / \rho )^{\al } \over (1 - \underline{p} + \underline{p} t / \rho )^{\ga }} \Big] _{1}^{\rho } / \int_{\rho }^{1} {1 \over \rho ^{\al }} {t^{\al - 1} \over (1 - \underline{p} + \underline{p} t / \rho )^{\ga + 1}} dt \non \\
&= - {1 - \underline{p} \over 1 - \overline{p}} \Big[ {t^{\al } \over \{ 1 - \overline{p} (1 - t) \} ^{\ga }} \Big] _{\rho }^{1} / \int_{\rho }^{1} {t^{\al - 1} \over \{ 1 - \overline{p} (1 - t) \} ^{\ga + 1}} dt \text{.} \non 
\end{align}
Then 
\begin{align}
\frac{ \int_{\underline{p}}^{\overline{p}} p^{\al - 1} (1 - p)^{\ga - \al - 1} dp }{ \int_{\underline{p}}^{\overline{p}} p^{\al - 1} (1 - p)^{\ga - \al } dp } &= {\ga \over \ga - \al } \Big[ 1 - {1 \over (1 - \overline{p} ) \ga } \Big[ {t^{\al } \over \{ 1 - \overline{p} (1 - t) \} ^{\ga }} \Big] _{\rho }^{1} / \int_{\rho }^{1} {t^{\al - 1} \over \{ 1 - \overline{p} (1 - t) \} ^{\ga + 1}} dt \Big] \non \\
&= {\ga \over \ga - \al } \Big[ 1 - {1 \over (1 - \overline{p} ) \ga } \Big[ {t^{\al } \over \{ 1 - \overline{p} (1 - t) \} ^{\ga }} \Big] _{\rho }^{1} / I( \al , \ga + 1, \underline{p} , \overline{p} ) \Big] \text{.} \non 
\end{align}
For part (iii), we have, by (\ref{lI_bothp1}), 
\begin{align}
\al \int_{\rho }^{1} {t^{\al - 1} \over (1 + \overline{r} t)^{\ga }} dt %
&= \Big[ {t^{\al } \over (1 + \overline{r} t)^{\ga }} \Big] _{\rho }^{1} + \ga \Big\{ \int_{\rho }^{1} {t^{\al - 1} \over (1 + \overline{r} t)^{\ga }} dt - \int_{\rho }^{1} {t^{\al - 1} \over (1 + \overline{r} t)^{\ga + 1}} dt \Big\} \non 
\end{align}
or 
\begin{align}
\ga \int_{\rho }^{1} {t^{\al - 1} \over (1 + \overline{r} t)^{\ga + 1}} dt &= \Big[ {t^{\al } \over (1 + \overline{r} t)^{\ga }} \Big] _{\rho }^{1} + ( \ga - \al ) \int_{\rho }^{1} {t^{\al - 1} \over (1 + \overline{r} t)^{\ga }} dt \text{.} \non 
\end{align}
Therefore, by the covariance inequality, 
\begin{align}
&\Big[ {t^{\al } \over (1 + \overline{r} t)^{\ga + 1}} \Big] _{\rho }^{1} / \int_{\rho }^{1} {t^{\al - 1} \over (1 + \overline{r} t)^{\ga + 1}} dt + \ga - \al + 1 \non \\
&= ( \ga + 1) \int_{\rho }^{1} {t^{\al - 1} \over (1 + \overline{r} t)^{\ga + 2}} dt / \int_{\rho }^{1} {t^{\al - 1} \over (1 + \overline{r} t)^{\ga + 1}} dt \non \\
&\ge ( \ga + 1) \int_{\rho }^{1} {t^{\al } \over (1 + \overline{r} t)^{\ga + 2}} dt / \int_{\rho }^{1} {t^{\al } \over (1 + \overline{r} t)^{\ga + 1}} dt \non \\
&= \Big[ {t^{\al + 1} \over (1 + \overline{r} t)^{\ga + 1}} \Big] _{\rho }^{1} / \int_{\rho }^{1} {t^{\al } \over (1 + \overline{r} t)^{\ga + 1}} dt + \ga - \al \non 
\end{align}
or 
\begin{align}
\Big[ {t^{\al + 1} \over \{ 1 - \overline{p} (1 - t) \} ^{\ga + 1}} \Big] _{\rho }^{1} / I( \al + 1, \ga + 1, \underline{p} , \overline{p} ) &\le 1 + \Big[ {t^{\al } \over \{ 1 - \overline{p} (1 - t) \} ^{\ga + 1}} \Big] _{\rho }^{1} / I( \al , \ga + 1, \underline{p} , \overline{p} ) \text{.} \non 
\end{align}
Furthermore, 
\begin{align}
\Big[ {t^{\al } \over \{ 1 - \overline{p} (1 - t) \} ^{\ga }} \Big] _{\rho }^{1} &< \Big[ {t^{\al + 1} \over \{ 1 - \overline{p} (1 - t) \} ^{\ga + 1}} \Big] _{\rho }^{1} \non 
\end{align}
and 
\begin{align}
\Big[ {t^{\al } \over \{ 1 - \overline{p} (1 - t) \} ^{\ga + 1}} \Big] _{\rho }^{1} &< \Big[ {t^{\al } \over \{ 1 - \overline{p} (1 - t) \} ^{\ga }} \Big] _{\rho }^{1} \text{.} \non 
\end{align}
Thus, 
\begin{align}
\Big[ {t^{\al } \over \{ 1 - \overline{p} (1 - t) \} ^{\ga }} \Big] _{\rho }^{1} / I( \al + 1, \ga + 1, \underline{p} , \overline{p} ) &< \Big[ {t^{\al + 1} \over \{ 1 - \overline{p} (1 - t) \} ^{\ga + 1}} \Big] _{\rho }^{1} / I( \al + 1, \ga + 1, \underline{p} , \overline{p} ) \non \\
&\le 1 + \Big[ {t^{\al } \over \{ 1 - \overline{p} (1 - t) \} ^{\ga + 1}} \Big] _{\rho }^{1} / I( \al , \ga + 1, \underline{p} , \overline{p} ) \non \\
&< 1 + \Big[ {t^{\al } \over \{ 1 - \overline{p} (1 - t) \} ^{\ga }} \Big] _{\rho }^{1} / I( \al , \ga + 1, \underline{p} , \overline{p} ) \text{.} \non 
\end{align}
This completes the proof. 
\hfill$\Box$

\bigskip

\noindent
{\bf Proof of Lemma \ref{lem:partialp}.} \ \ We have 
\begin{align}
{\pd \over \pd p} E_{p}^{X} [ \varphi (X) ] %
&= E_{p}^{X} \Big[ \Big( {X \over p} - {n - X \over 1 - p} \Big) \varphi (X) \Big] = {1 \over p (1 - p)} E_{p}^{X} [ (X - n p) \varphi (X) ] \non \\
&= {1 \over p} E_{p}^{X} [ X \{ \varphi (X) - \varphi (X - 1) \} ] \text{,} \non 
\end{align}
where the third equality follows from the lemma of Johnson (1987). 
Therefore, 
\begin{align}
\Big( {\pd \over \pd p} \Big) ^2 \{ p E_{p}^{X} [ \varphi (X) ] \} &= {\pd \over \pd p} [ E_{p}^{X} [ \varphi (X) ] + E_{p}^{X} [ X \{ \varphi (X) - \varphi (X - 1) \} ] ] \non \\
&= {\pd \over \pd p} E_{p}^{X} [ (X + 1) \varphi (X) - X \varphi (X - 1) ] \non \\
&= {1 \over p} E_{p}^{X} [ X \{ (X + 1) \varphi (X) - 2 X \varphi (X - 1) + (X - 1) \varphi (X - 2) \} ] \text{,} \non 
\end{align}
which is the desired result. 
\hfill$\Box$

\bigskip

\noindent
{\bf Proof of Lemma \ref{lem:log_Jensen}.} \ \ We have 
\begin{align}
E^T [ \log (1 - T) ] &= - \sum_{k = 1}^{\infty } {1 \over k} E^T [ T^k ] = - \mu - {\si ^2 + \mu ^2 \over 2} - \sum_{k = 3}^{\infty } {1 \over k} E^T [ T^k ] \non \\
&\le - \mu - {\si ^2 + \mu ^2 \over 2} - \sum_{k = 3}^{\infty } {1 \over k} \mu ^k = \log (1 - \mu ) - {\si ^2 \over 2} \text{,} \non 
\end{align}
where the inequality follows from Jensen's inequality. 
\hfill$\Box$

\section*{Acknowledgments}
I would like to thank Professor Tatsuya Kubokawa for his encouragement. 
In particular, Section \ref{sec:relation_Po} is based on his comments. 
Research of the author was supported in part by JSPS KAKENHI Grant Number JP20J10427 from Japan
Society for the Promotion of Science. 

\end{document}